\documentclass[epic,11pt]{article}
\usepackage{amsfonts,amsmath,amscd,amssymb,amsthm}
\usepackage{fullpage}
\usepackage[pdftex]{pict2e}

\usepackage[alphabetic,initials]{amsrefs}
\usepackage[all]{xy}
\usepackage{hyperref} 
\usepackage{mathrsfs}
\input xy%
\usepackage{fullpage}
\usepackage[pdftex]{pict2e}
\def\0{{\bar 0}}
\def\1{{\bar 1}}

\def\Z{{\mathbb Z}}
\def\X{{\mathbb X}}
\def\Y{{\mathbb Y}}
\def\F{{\mathbb F}}

\def\B{{\mathbb B}}
\def\G{{\mathbb G}}

\def\cpa{{\operatorname{c.empty}}} 
\def\cda{{\operatorname{c.full}}} 
\def\rpa{{\operatorname{r.empty}}} 
\def\rda{{\operatorname{r.full}}}

\def\H{{\mathbb H}}

\def\W{\mathtt{W}}

\def\ur{{\operatorname{r  }}}
\def\ud{{\operatorname{d}}}

\def\span{\operatorname{span\;}}

\def\str{{\operatorname{red}}}
 
\def\trans{{\operatorname{transpose}}} 
\def\dist{{\operatorname{dist}}}

\def\iso{{\operatorname{iso}}}

\newcommand{\ap}{\approx}

\newcommand{\ttk}{\mathtt{k}}
\newcommand{\ttr}{\mathtt{r}}
\newcommand{\ttw}{\mathtt{w}}
\newcommand{\ttc}{\mathtt{c}}

\newcommand{\ttd}{\mathtt{d}}

\newcommand{\ghs}{{\gth(\gs)}}

\newcommand{\diag}{diag}

\newcommand{\Lgo}{L(\stackrel{{\rm _o}}{\fg})}
\newcommand{\Lgh}{\widehat{L}(\stackrel{{\rm _o}}{\fg})}
\newcommand{\itemi}{\item[{{\rm(i)}}]}
\newcommand{\itemii}{\item[{{\rm(ii)}}]}

\newcommand{\itema}{\item[{{\rm$($a$)$}}]}
\newcommand{\itemb}{\item[{{\rm$($b$)$}}]}
\newcommand{\itemc}{\item[{{\rm$($c$)$}}]}
\newcommand{\itemd}{\item[{{\rm$($d$)$}}]}
\newcommand{\iteme}{\item[{{\rm$($e$)$}}]}

\newcommand{\itemo}{\item[{}]}
\newcommand{\noi}{\noindent}

\newcommand{\ga}{\alpha}
\newcommand{\gb}{\beta}
\newcommand{\gc}{\gamma}

\newcommand{\Gc}{\Gamma}
\newcommand{\Gl}{\Lambda}
\newcommand{\Gd}{\Delta}

\newcommand{\gd}{\delta}

\newcommand{\gs}{\sigma}

\newcommand{\go}{\omega}

\newcommand{\gt}{\tau}
\newcommand{\gz}{\zeta}

\newcommand{\gl}{\lambda}

\newcommand{\gr}{\rho}

\newcommand{\gep}{\epsilon}
\newcommand{\gth}{\theta}
\newcommand{\op}{\oplus}

\newcommand{\mr}{{-_\ttr}}
\newcommand{\pr}{{+_\ttr}}

\newcommand{\omc}{{-_\ttc}}
\newcommand{\opc}{{+_\ttc}}

\newcommand{\fg}{\mathfrak{g}}\newcommand{\fgl}{\mathfrak{gl}}
\newcommand{\fsl}{\mathfrak{sl}}

\newcommand{\fh}{\mathfrak{h}}

\newcommand{\fb}{\mathfrak{b}}

\newcommand{\fB}{\mathfrak{B}}

\newcommand{\fS}{\mathfrak S}
\newcommand{\fT}{\mathfrak T}

\newcommand{\ff}{\footnote}
\newfont{\eufm}{eufm10 scaled\magstep1}

 \newcommand{\ti}{\times}

\newcommand{\bcu}{\bigcup}

\newcommand{\cB}{\mathcal{B}}

\newcommand{\cE}{\mathcal{E}}

\newcommand{\cG}{\mathcal{G}}

\newcommand{\cS}{\mathcal{S}}

\newcommand{\cX}{\mathcal{X}}

\newcommand{\cW}{\mathfrak{W}}

\newcommand{\ttx}{\mathtt{x}}

\newcommand{\Sh}{{\bf{Sh}}}
\newcommand{\sh}{{\bf{sh}}}
\newcommand{\ey}{\end{eqnarray}}
\newcommand{\by}{\begin{eqnarray}}
\newcommand{\nn}{\nonumber}

\newcommand{\bco}{\begin{conjecture}}
\newcommand{\ba}{\begin{alg}}
\newcommand{\ea}{\end{alg}}
\newcommand{\eco}{\end{conjecture}}
\newcommand{\bpf}{\begin{proof}}
\newcommand{\epf}{\end{proof}}
\newcommand{\bt}{\begin{theorem}}

\newcommand{\et}{\end{theorem}}

\newcommand{\br}{\begin{rem}}
\newcommand{\er}{\end{rem}}
\newcommand{\brs}{\begin{rems}}
\newcommand{\ers}{\end{rems}}
\newcommand{\bi}{\begin{itemize}}
\newcommand{\ei}{\end{itemize}}
\newcommand{\bl}{\begin{lemma}}
\newcommand{\bsul}{\begin{sublemma}}
\newcommand{\esul}{\end{sublemma}}
\newcommand{\bp}{\begin{proposition}}
\newcommand{\be}{\begin{equation}}
\newcommand{\bc}{\begin{corollary}}
\newcommand{\bexs}{\begin{examples}}
\newcommand{\eexs}{\end{examples}}
\newcommand{\bexa}{\begin{example}}
\newcommand{\eexa}{\end{example}}
\newcommand{\bex}{\begin{exercise}}
\newcommand{\eex}{\end{exercise}}
\newcommand{\btab}{\begin{tab}}
\newcommand{\etab}{\end{tab}}
\newcommand{\el}{\end{lemma}}
\newcommand{\ep}{\end{proposition}}
\newcommand{\ee}{\end{equation}}
\newcommand{\ec}{\end{corollary}}
\newcommand{\Bc}{\begin{center}}
\newcommand{\Ec}{\end{center}}
\newcommand{\bh}{\begin{hyp}}
\newcommand{\eh}{\end{hyp}}
\newcommand{\bhs}{\begin{hyps}}
\newcommand{\ehs}{\end{hyps}}
\newcommand{\bd}{\begin{dfn}}
\newcommand{\ed}{\end{dfn}}

\newcommand{\bn}{\begin{notn}}
\newcommand{\en}{\end{notn}}
\newcommand{\ogd}{\overline{\delta}}

\newcommand{\fgm}{\stackrel{{\rm o}}{\fg}^{-\gth}}

\newcommand{\fgo}{\stackrel{{\rm o}}{\fg}}
\newcommand{\fgu}{\underline{\fgo}}
\newcommand{\pou}{\underline{\pigo}}
\newcommand{\pvu}{\underline{\pigo}^\vee}
\newcommand{\fag}{\underline{A}^0}
\newcommand{\fto}{\underline{\gt}^0}
\newcommand{\Hbo}{\stackrel{{\rm o}}{\cB}}
\newcommand{\pigo}{\stackrel{{\rm o}}{\Pi}}
\newcommand{\fbo}{\stackrel{\rm o}{\fb}}{}
{}

\newcommand{\pigov}{\stackrel{\rm o}{\Pi^\vee}}

\newcommand{\fgt}{\stackrel{{\rm o}}{\fg}^\gth}

\newcommand{\fho}{\stackrel{\rm o}{\fh}}
\newcommand{\fhs}{\stackrel{{\rm _o}}{\fh}^*}

\usepackage[aligntableaux=center]{ytableau}
\usepackage{tabularray}

\begin{document}
\title{Table of Contents}

\newtheorem*{bend}{Dangerous Bend}

\newtheorem{thm}{Theorem}[section]
\newtheorem{hyp}[thm]{Hypothesis}
 \newtheorem{hyps}[thm]{Hypotheses}
\newtheorem{notn}[thm]{Notation}

  \newtheorem{rems}[thm]{Remarks}

\newtheorem{conjecture}[thm]{Conjecture}
\newtheorem{theorem}[thm]{Theorem}
\newtheorem{theorem a}[thm]{Theorem A}
\newtheorem{example}[thm]{Example}
\newtheorem{examples}[thm]{Examples}
\newtheorem{corollary}[thm]{Corollary}
\newtheorem{rem}[thm]{Remark}
\newtheorem{lemma}[thm]{Lemma}
\newtheorem{sublemma}[thm]{Sublemma}
\newtheorem{cor}[thm]{Corollary}
\newtheorem{proposition}[thm]{Proposition}
\newtheorem{exs}[thm]{Examples}
\newtheorem{ex}[thm]{Example}
\newtheorem{exercise}[thm]{Exercise}
\numberwithin{equation}{section}%
\setcounter{part}{0}
\newcommand{\drar}{\rightarrow}
\newcommand{\lra}{\longrightarrow}
\newcommand{\rra}{\longleftarrow}
\newcommand{\dra}{\Rightarrow}
\newcommand{\dla}{\Leftarrow}
\newcommand{\rl}{\longleftrightarrow}

\newtheorem{Thm}{Main Theorem}


\newtheorem*{thm*}{Theorem}
\newtheorem{lem}[thm]{Lemma}
\newtheorem*{lem*}{Lemma}
\newtheorem*{prop*}{Proposition}
\newtheorem*{cor*}{Corollary}
\newtheorem{dfn}[thm]{Definition}
\newtheorem*{defn*}{Definition}
\newtheorem{notadefn}[thm]{Notation and Definition}
\newtheorem*{notadefn*}{Notation and Definition}
\newtheorem{nota}[thm]{Notation}
\newtheorem*{nota*}{Notation}
\newtheorem{note}[thm]{Remark}
\newtheorem*{note*}{Remark}
\newtheorem*{notes*}{Remarks}
\newtheorem{hypo}[thm]{Hypothesis}
\newtheorem*{ex*}{Example}
\newtheorem{prob}[thm]{Problems}
\newtheorem{conj}[thm]{Conjecture}

\title{Young diagrams, Borel subalgebras and Cayley graphs}
\author{Ian M. Musson
 \\Department of Mathematical Sciences\\
University of Wisconsin-Milwaukee\\ email: {\tt
musson@uwm.edu}}
\maketitle
\begin{abstract} Let $\ttk$ be an algebraically closed field of characteristic zero and $n, m$ coprime positive integers.  Let ${\stackrel{{\rm o}}{\mathfrak{g}}}$ be the Lie superalgebra $\fsl(n|m)$ and let  
$\mathfrak T_{iso}$ 
be the groupoid 
introduced by 
 Sergeev and Veselov \cite{SV2} with base the set of odd roots of ${\stackrel{{\rm o}}{\mathfrak{g}}}$.  
We show the Cayley graphs for three actions of 
$\mathfrak T_{iso}$ are isomorphic,  
These actions originate in quite different ways.  Consider the set $X$ of Young diagrams contained in a rectangle with $n$ rows and $m$ columns. By adding or deleting rows and columns from certain diagrams and keeping track of the total number of 
boxes added or deleted, we obtain an equivalence relation on $X\times  \Z$ such that   $\mathfrak  T_{iso}$ acts on the  set of equivalence classes $[X\times  \Z]$. 
We compare the action on $[X\times  \Z]$ to an action on Borel subalgebras of  the affinization $\Lgh$ of  ${\stackrel{{\rm o}}{\mathfrak{g}}}$ which are related by odd reflections. The third action 
comes from an action of $\mathfrak T_{iso}$ on $\ttk^{n|m}$  
defined  by 
 Sergeev and Veselov,  motivated by deformed quantum Calogero-Moser problems \cite{SV1}.  This action will be considered in \cite{M24}.
\end{abstract}




\section{Introduction}\label{itr} 
Since their introduction nearly a century ago \cite{BT}, groupoids have become an indispensable tool in Mathematics and Physics. For concise and readable surveys  see \cite{BR}, \cite{We}.   A famous example is the van Kampen Theorem in  topology,  \cite{CF}, \cite{Hi}.
A widely held point of view is that, in certain situations groupoids can be  more effective than groups in the study of symmetry. A key concept in this respect is the action of  a groupoid  $\cG$ on a set  $\X$. 
For a set $\X$, the{\it  symmetric groupoid} $\cS(\X)$  on  $\X$ has base consisting of all subsets of $\X$ and morphisms all bijections between subsets.
The groupoid $\cG$ {\it  acts on} $\X$ if there is a functor $F:\cG\lra\cS(\X)$.  If $X$ is an affine algebraic variety, we can study the invariant algebra consisting of all polynomial functions on $X$ that are constant on $\cG$-orbits.  Actions and  invariants of groupoids are far less understood than their analogs for groups.  Interesting examples arise in the representation theory of a basic classical Lie superalgebra $\fg$. 
 Indeed 
 Sergeev and Veselov \cite{SV2} introduced an action of the Weyl groupoid 
$\cW$, whose invariants are closely related to the Grothendieck group of finite dimensinonal graded $\fg$-modules.   The center $Z(\fg)$ of $U(\fg) $ arises as the invariant ring for a related action of $\cW$, see 
\cite{M22}.  
The Weyl groupoid is constructed using the usual Weyl group and a certain subgroupoid
$\mathfrak T_{iso}$.  

\subsection{ Groupoids and Cayley graphs }\label{itr1} 
Throughout  $[k]$ denotes  the set of the first $k$ positive integers and iff means if and only if.  We introduce a tool that may be used to compare two such actions of the same groupoid. 
Suppose $S$ is a set of morphisms in $\cG$. The {\it Cayley graph} $\Gc(S,\cG)$ of $F$ is the 
directed graph with vertex set $\X$ and an edge $x\lra F(s)x$ whenever $s\in  S$ and $F(s)x$ is defined.  
If $ \cG$ is a group acting on itself by the regular representation, we get the usual Cayley graph.
\\ \\We assume $\cG=\mathfrak  T_{iso}$ 
and take $S$ to be all morphisms in $\cG$.  If $x\in \X$, we refer to the (vertices in the) connected component of the Cayley graph as the {\it orbit} of $x$.  
Thus what is important to us is the functor $F$. 
  Therefore we  denote $\Gc(S,\cG)$ instead by  
$\Gc(F)$.
\\ \\
If $\go:\X\lra \Y$ is a bijection of sets, there is an induced isomorphism $\cS(\X)\lra \cS(\Y)$  also denoted $\go$.  The morphism $g:S\lra T$ in  $\cS(\X)$ corresponds to $\go(g)= \go g \go^{-1}:\go(S)\lra \go(T)$.
   Now suppose we have a commutative diagram of groupoids and functors
\be \label{bby}
\xymatrix@C=2pc@R=1pc{
\cG\ar@{<->}_= [dd]&&
\cS(\X)\ar@{<-}_F[ll]: \ar@{->}^{\go}[dd]&\\ \\
\cG \ar@{-}[rr]_B &&
\cS(\Y)\ar@{<-}[ll] &}
\ee
Then for any $x\in \X$ and morphism $g$ in $\cG$
$$\go(F(g)x)=B(g)\go(x).$$  Thus $\go$ induces an isomorphism of Cayley graphs 
$$\Gc(F) \lra \Gc(B).$$
We say that $F$ acts {\it transitively} if $\X$  is a single $\cG$-orbit, compare \cite{BR} Section 4. If $F$ and $B$ are transitive in \eqref{bby}, we say that $\go$ is a $\cG$-{\it equivariant bijection from 
$\X$ to $\Y$}. 
The aim of this paper and its sequel \cite{M24} is to construct various commutative diagrams as in \eqref{bby}.

\subsection{The Weyl groupoid}\label{aog}
The Lie superalgebras $\fgl(n|m)$ and  $\fsl(n|m)$ have the same root system $\Gd$
and there is an obvious bijection between Borel subalgebras of  $\fgl(n|m)$  and $\fsl(n|m)$.  
It will be convenient to work with
$\fgo= \fgl(n|m)$ except in  Section \ref{akmls}
where $\fgo=\fsl(n|m)$. 
We consider the Weyl groupoid $\cW$ 
 defined in \cite{SV2} 
using the root system  $\Gd$. Let $\mathfrak T_{iso}$ be the groupoid with base 
$\Gd_{iso}$ the set of all the isotropic roots in $\Gd.$  
The non-identity  morphisms  are  $\gr_{\alpha}:\alpha \rightarrow -\ga$, $\ga\in\Gd_{iso}$. The Weyl group $W$ acts on $\mathfrak T_{iso}$ in a natural way: $\alpha \rightarrow w(\alpha),\,
\gr_{\alpha} \rightarrow \gr_{w(\alpha)}$.  The {\it Weyl groupoid} $
\cW$ is defined
by 
\be \label{fgh}  \mathfrak{W} = W \coprod W \ltimes \mathfrak T_{iso},\ee     the disjoint union of the group $W$ considered as a groupoid with a single point base $[W]$ and the semidirect product $W \ltimes \mathfrak T_{iso}$. \ff{According to Brown \cite{BR}, semidirect products of groupoids were first considered by Ehresmann  \cite{E}.}
We work almost exclusively with the groupoid $\mathfrak T_{iso}$.
If $\ga\in \Gd^+_1$ we call $\gr_\ga$ a {\it positive morphism} and likewise for the morphisms $t_\ga =  F(\gr_\ga)$, $p_\ga =  G(\gr_\ga)$ and $r_\ga =  B(\gr_\ga)$ 
 defined in Subsection \ref {iog}.

\subsection{ A first example} \label{itr2} 
A {\it partition} $\gl$ into $n$ {\it parts} is a  descending sequence 
\be \label{pdef}  
\gl=(\gl_1, \gl_2, \ldots, \gl_n)\ee of non-negative integers.  
We call $\gl_i$ the $i^{th}$ {\it part} $\gl$.  
A partition is identified with the Young diagram it defines, but we say $\gl$ {\it is a partition} or $\gl$  {\it is a  Young diagram} according to which attributes of $\gl$ we wish to emphasize. 
Given positive integers $m, n$, 
let $X$ be the set of Young diagrams that fit inside (the bottom left corner of) 
a  rectangle {\bf R}  with $n$ rows and $m$ columns.  Boxes in 
{\bf R} correspond to root spaces in  the 
upper triangular block $\fg^+_1$ of $\fgl(n|m)$.  Accordingly the rows of   
{\bf R} are indexed from the top down by $\gep_1,\ldots,\gep_n$  and columns of {\bf R} are indexed  from  left to right by  $\gd_1,\ldots,\gd_m$.  Thus the box in row $\gep_i$ and column $\gd_j$
corresponds to the root space $\fg^{\gep_i-\gd_j}$ and we call this box $\fB(\gep_i-\gd_j)$, or often just $\gep_i-\gd_j$.  The set of roots of 
$\fg^+_1$ is $\Gd_1^+ = \{\gep_i-\gd_j|i\in [n], j \in [m]\}.$
\\ \\
We also write 
\be \label{prf}\gl=(m^{a_m},  \ldots, i^{a_i}, \ldots)\ee
 to indicate that the part $i$ occurs with multiplicity ${a_i}$ in $\gl$.   Set $|\gl| = \sum_{i=1}^n \gl_i$. 
For a partition $\gl$,  the {\it dual partition} $\gl'$  has $j^{th}$ part 
$\gl'_j = |\{i|\gl_i \ge j\}|$. 
\\ \\ 
 Also let ${\stackrel{{\rm o}}{\cB}}$ be the set of Borel subalgebras $\fbo$ of $\fgl(n|m)$ such that 
$\fbo_0$ is the Borel subalgebra of $\fg$ consisting of two diagonal blocks of 
upper triangular matrices.  There is a natural bijection $\go:X\rl {\stackrel{{\rm o}}{\cB}}$
which is  $\mathfrak T_{iso}$-equivariant for the folllowing actions.  
\bi
\item
$\mathfrak T_{iso}$  acts  on $X$ by 
adding and deleting corners to and from a Young diagram. 
\item 
$\mathfrak T_{iso}$  acts on  ${\stackrel{{\rm o}}{\cB}}$ by odd reflections.
\ei
Denote the empty Young diagram $(0^n)$ by $\emptyset$ and the distinguished Borel subalgebra of $\fgo$ by ${\stackrel{{\rm o}}{\fb}}{}^\dist$. We require that 
 $\go(\emptyset)={\stackrel{{\rm o}}{\fb}}{}^\dist$.
 For details 
see Subsection \ref{iog}. 
Let
\be \label{fnx} F:\mathfrak T_{iso}\lra\cS({X}), \quad  B:\mathfrak  T_{iso} \lra \cS( {\stackrel{{\rm o}}{\cB}})\ee
be the functors corresponding to the above actions. 

\subsection{ Extending  the action of $\mathfrak T_{iso}$} 
\label{itr3}
We say that $\gl \in X$ is {\it row empty} if the highest row of $\gl$ is empty 
 and {\it row full} if the lowest row of $\gl$ contains the maximum possible 
number $m$ of boxes. The terms {\it column empty} 
 and {\it column full} are defined similarly. More formally \be\label{yyt}X(\rda) = \{\gl\in X|\gl_1=m\}, \quad X(\rpa) = \{\gl\in X|\gl_n=0\}.\ee
\be\label{yys}X(\cda) = \{\gl\in X|\gl'_1=n\}, \quad X(\cpa) = \{\gl\in X|\gl'_m=0\}.\ee
If $\gl \in X(\rda)$, 
let   $\gl^\mr =  (\gl_2, \ldots,\gl_n,0 )$  be the partition obtained from removing the first part from $\gl$  
and if $\gl \in X(\rpa)$, set 
$\gl^\pr =(m,\gl_1,\gl_2,\ldots ,\gl_{n-1}).$\ff{For the notation, might help to think that if something is empty, then we can add to it.}
As a Young diagram 
${\gl}^\mr $
is obtained from  $\gl$ by deleting the bottom row.   We have  inverse bijections
$X(\rda) \rl X(\rpa)$ given by  $\gl \lra \gl^\mr$ and 
$\mu^\pr \rra \mu.$\\ \\
Similarly, if $\gl \in X(\cda)$, 
let   $\gl^\omc =  (\gl_1 -1,\gl_2-1, \ldots,\gl_n -1)$  
and if $\gl \in X(\cpa)$, set 
$\gl^\opc =(\gl_1 +1,\gl_2+1, \ldots,\gl_n +1)).$  
As a Young diagram 
${\gl}^\omc $
is obtained from  $\gl$ by deleting the first  column.   We have  inverse bijections
$X(\cda) \rl X(\cpa)$ given by  $\gl \lra \gl^\omc$ and 
$\mu^\opc \rra \mu.$
\noi Let $\sim$ be the smallest equivalence relation on 
$X\ti \Z$ such that $(\gl,k)\sim ({\gl}^\pr,k-m) $ if $
 \gl\in X(\rpa)$ and $(\gl,k)\sim ({\gl}^\opc,k-n) $
{ if }
 $\gl\in X(\cpa)$. 
Equivalently
\be \label{1kr}(\gl,k)\sim ({\gl}^\mr,k+m) \mbox{ if } \gl\in X(\rda)\ee and 
\be \label{2kr}(\gl,k)\sim ({\gl}^\omc,k+n) \mbox{ if } \gl\in X(\cda).\ee
Denote the equivalence class of $(\gl,k)\in X\ti Z$ under $\sim$ by $[\gl,k]$ and let  
$[X\ti Z]$ be the set of equivalence classes.  In the pair $({\gl},k)  \in X\ti \Z$, 
 we call $k$ the {\it rotation number}.
There is a compatible 
equivalence relation on 
${\stackrel{{\rm o}}{\cB}}\ti \Z.$  We show in Section \ref{APG} that the action of $\mathfrak T_{iso}$  can be extended to the set of equivalence classes  $[X \ti \Z]$ and   $[{\stackrel{{\rm o}}{\cB}} \ti \Z]$. 
Using the functors $F, B$ from 
\eqref{fnx}, we construct functors
\be \label{gnx} \F:\mathfrak T_{iso}\lra\cS[X \ti \Z], \quad  \B:\mathfrak  T_{iso} \lra \cS[{\stackrel{{\rm o}}{\cB}} \ti \Z]
\ee
and  in Corollary \ref{sz} we show there is a 
$\mathfrak T_{iso}$-equivariant bijection 
\be \label{fx}[X \ti \Z] \lra [{\stackrel{{\rm o}}{\cB}} \ti \Z].\ee The equivalence classes $[\emptyset, 0]$ and $[{\stackrel{{\rm o}}{\fb}}{}^\dist,0]$ correspond  under this bijection. 
\subsection{ Borel subalgebras of the affinization $\Lgh$} \label{itr5} 
The affinization  $\Lgh$ of $\fgo$  has infinitely many  Borel subalgebras, which arise in the following way.
The process of affinization produces subalgebras ${\stackrel{\rm o}{\fg}}(k)$ for $k\in\Z,$ all  isomorphic to $\fgo$.  In addition for each Young diagram $\gl$ there is a Borel subalgebra ${\stackrel{{\rm o}}{\fb}}{}(\gl,k)$ 
of ${\stackrel{\rm o}{\fg}}(k)$ as in \cite{M101}  Chapter 3 and each subalgebra ${\stackrel{{\rm o}}{\fb}}{}(\gl,k)$ extends to a unique Borel subalgebra of $\Lgh$. 
Also by affinization any Dynkin-Kac diagram for ${\stackrel{\rm o}{\fg}}(k)$ is  embedded in a corresponding diagram for  $\Lgh$ by adding an 
extending node.   Under the conditions stated below, by deleting a different node we
obtain 
\by \label{f2}   \mbox{ a Borel subalgebra }
&{\stackrel{{\rm o}}{\fb}}{}(\gl^\mr,k+m)& \mbox{ of } \quad {\stackrel{\rm o}{\fg}}(k+m) \quad \mbox{ if }\quad  \gl\in X(\rda)\\\nn \mbox{ a Borel subalgebra }
&{\stackrel{{\rm o}}{\fb}}{}(\gl^\pr,k-m)& \mbox{ of } \quad {\stackrel{\rm o}{\fg}}(k-m) \quad \mbox{ if } \quad \gl\in X(\rpa)\\ \nn
\mbox{ a Borel subalgebra }
&{\stackrel{{\rm o}}{\fb}}{}(\gl^\omc,k+n)& \mbox{ of } \quad {\stackrel{\rm o}{\fg}}(k+n) \quad \mbox{ if } \quad \gl\in X(\cda)\\\nn \mbox{ a Borel subalgebra }
&{\stackrel{{\rm o}}{\fb}}{}(\gl^\opc,k-n)& \mbox{ of } \quad  {\stackrel{\rm o}{\fg}}(k-n) \quad \mbox{ if } \quad \gl\in X(\cpa)\\\nn
\ey  For details see Subsection \ref{afz}. 
If $n,m$ are coprime, then starting from the distinguished Borel 
${\stackrel{{\rm o}}{\fb}}{}(\emptyset, 0)$ of ${\stackrel{\rm o}{\fg}}(0)=\fgo$, and adding/deleting rows, columns and corners we can reach the Borel  subalgebra 
${\stackrel{{\rm o}}{\fb}}{}(\mu,\ell)$ of $\fgo(\ell)$ for any pair $(\mu,\ell) \in {\stackrel{{\rm o}}{\cB}}\ti \Z.$
Let $\cB$ be the set of Borel subalgebras of  $\Lgh$ that are extensions of  Borels in the subalgebras ${\stackrel{\rm o}{\fg}}(k)$ for  $k\in \Z.$  
There is a functor 

\be \label{Fnx}{\bf B}:\mathfrak  T_{iso}\lra \cS(\cB).\ee 
This action is transitive: If $\fb$ is the Borel subalgebra of $\Lgh$ extended from the distinguished subalgebra of ${\stackrel{\rm o}{\fg}}(0)$, then $\cB$ is  the  $\mathfrak  T_{iso}$-orbit  
of 
$\fb$.

	
\subsection{ Main result} \label{itr6} 
We consider the following functors and orbits \bi \item The functor $\F:\mathfrak T_{iso}\lra\cS[X \ti \Z]$ from \eqref{gnx}
\item 
Let 
${\bf B}:\mathfrak  T_{iso} \lra \cS( \cB)$  be the functor 
defined in \eqref{Fnx}.
\item Let $SV:\mathfrak  T_{iso} \lra \cS(\Z^{n|m})$  and let {\bf O}
be the  $\mathfrak  T_{iso}$-orbit  
of 
$\Gl_0 \in \Z^{n|m}$ as defined in \cite{M24}, see also \cite{SV101} Equation (14).
\ei
Our main result relates the first two of these functors. 
\bt \label{iir}  Assume  $m, n$ are coprime and  $\fgo = \fsl(m|n)$ there is a 
$\mathfrak T_{iso}$-equivariant bijection $[X \ti \Z]\lra\cB$
such that $[\emptyset, 0]$ corresponds to $\fb$. 
\et \noi  
An important point is that the 
the Borels ${\stackrel{{\rm o}}{\fb}}{}(\gl,k)$ and ${\stackrel{{\rm o}}{\fb}}{}(\gl^\mr,k+m)$ have the same extension to $\Lgh$ if $\gl\in X(\rda)$.  
An analogous statement holds if $\gl\in X(\cda)$.  
This leads in 
Theorem \ref{srb}, to a 
$\mathfrak T_{iso}$-equivariant bijection  
$\cB\rl [{\stackrel{{\rm o}}{\cB}} \ti \Z]$.  Theorem \ref{iir} follows from this and 
\eqref{fx}. 

\subsection{ The action of Sergeev and Veselov} 
\label{itr4} 
The Lie superalgebras $\fgl(n|m)$ and  $\fsl(n|m)$ have the same root system $\Gd$
and there is an obvious bijection between Borel subalgebras of  $\fgl(n|m)$  and $\fsl(n|m)$.  
It will be convenient to work with
$\fgo= \fgl(n|m)$ except in Section \ref{akmls}
where $\fgo=\fsl(n|m)$.   
In  \cite{SV101} Sergeev and Veselov defined 
an action of the   $\mathfrak{W}$  on $\ttk^{n|m}$.  In  \cite{M24} we reformulate their definition so that 
 $\mathfrak{W}$  acts on $\Z^{n|m}$ via a functor
 $$SV:\mathfrak {W}\lra \cS(\Z^{n|m}).$$ The restriction of $SV$
to $\mathfrak  T_{iso}$ is also denoted $SV$. 
Regarding the third functor,
in  \cite{M24} we show
\bt \label{ipr}  If $m, n$ are coprime,  there is a 
$\mathfrak T_{iso}$-equivariant bijection $x:[X \ti \Z] \lra{\bf O}$
such that $x(\emptyset, 0) =   \Gl_0$. 
\et \noi  
This paper is organized as follows. Section \ref{fc} contains preliminary results on Young diagrams, groupoids and Borel subalgebras. 
In Section \ref{APG}   we explain how the action of $\mathfrak  T_{iso}$ can be extended to an action on a set of equivalence classes. 
In Section \ref{akmls},  we discuss 
Borel subalgebras of 
the affinization $\Lgh$. 
In Examples \ref{iob}, \ref{E3} and  \ref{E4} 
we illustrate many of our results in the  case
where $(n,m)=(2,3)$.
\\ \\
We study several objects which are very closely related: Young diagrams, partitions, words, shuffles and Borel subalgebras. We often assume that $w, \gl$ and $\gs$ are related as in Equation \eqref{bji} and then use whichever notation is most convenient.   For example the same Borel subalgebra of $\fgo$ may be denoted by $\fbo(\gs)$ or $\fbo(\gl)$.  The alternative to this abuse of notation seems to be to introduce more notation.  
\\ \\
Groupoids are used to study Borel subalgebras of affine Lie superalgebras in 
\cite{GHS}, \cite{GHS2} and \cite{GK}.  The  motivation and methods in these papers are quite different from ours. I thank Lukas Bonfert and  Maria Gorelik for some helpful correspondence.

\section{Young Diagrams, Groupoids and Borel subalgebras of $\fgl(n|m)$}\label{fc}
\subsection{The Lie superalgebra $\fgl(n|m)$}\label{Lsg}
Let $\fgo=\fgl(n|m)$.  The Cartan subalgebra $\fho$ of diagonal matrices in $\fgo$ is identified with $\ttk^{n|m}$.  
For $i\in [m+n], $  define $\gep_i\in \fhs$ so that $\gep_i(h) = h_i$ where  $h=\diag(h_1,\ldots,h_{m+n})\in \;\fho$.  For $i\in [m]$, set $i'=i+n$ and  $\gd_{i} =\gep_{i'}$. 
We have the usual permutation action  of $W$ on $\ttk^{n|m}$.  
The bilinear form $(\;,\;)$ is defined on $\fh^*$ is defined on the basis 
$\epsilon_{1},\ldots,   \epsilon_{n}, \delta_{1}, \ldots,\delta_{m}  $ by
\be \label{edform}(\epsilon_i,\epsilon_j) = \delta_{i,j} = - (\delta_i,\delta_j)\nn\ee 
and $(\epsilon_{i},\gd_{j})=0$ for all relevant indices $i,j$.
We write  $\Gl \in \ttk^{n|m}$ in the form
\be \label{ezr} \Gl = 
(a_1, \ldots ,a_n|b_1,\ldots,b_{m}) =  \sum_{i=1}^n  a_i
\gep_i - \sum_{j=1}^m  b_j \gd_j .\ee
If $\Gl$ is as above we have 
\be \label{exy} (\Gl,\gep_i- \gd_j)=  a_i-b_j. \ee 

\subsection{Rotation operators}\label{roo}
Define the permutation  $\nu$ of $\{1, \ldots, n\}$ by 
$\nu(k) = k-1 \mbox{ mod } n$. 
We extend $\nu$ to  an operator on $\ttk^{n|m}$  by setting
\be \label{e67l}\nu(a_1, \ldots, a_n|b_1, \ldots, b_m) = (a_{\nu(1)}, \ldots, a_{\nu(n)}|b_1, \ldots, b_m) \ee
Similarly define the permutation 
$\eta $ of $\{1', \ldots, m'\}$ by 
$\eta(k') = (k+1)' \mbox{ mod } m$ and set  
\be \label{e67l}\eta(a_1, \ldots, a_n|b_1, \ldots, b_m) = (a_{ 1}, \ldots, a_{n}|b_{\eta(1)}, \ldots, b_{\eta(m)}) \ee
Now \eqref{e67l} implies 
\be \label{e67m}\nu(\gep_i) = \gep_{i+1}  \mbox{ for } i\in [n], \quad \eta(\gd_{j+1}) = \gd_{j}  \mbox{ for } j\in [m].\ee
and 
$$\nu(\gd_j) = \gd_j ,\quad \eta(\gep_i)=\gep_{i} 
\mbox{ for } j\in [m], \quad 
 i\in [n].$$
We call $\nu$ and $\eta$ {\it rotation operators}. 

\subsection{Young diagrams and Borel subalgebras}\label{yd}  
 \subsubsection{Young diagrams, partitions and corners}\label{ypc}
In what follows we view a Young diagram $\gl$ as a set of boxes in {\bf R}.
If $\ga\in\Gd_1^+$,we say the box  $\fB=\fB( \ga)$ is an  {\it outer corner} of $\gl$ if $
  \fB \notin \gl$ 
and the addition of $ \fB$ to $\gl$ creates a new Young diagram.
Similarly box $ \fB$ is an  {\it inner corner} of $\gl$ if $
 \fB \in \gl$ 
and the removal of $ \fB$ from $\gl$ creates a new Young diagram.
Let 
$$X_\ga = \{\gl\in X| \fB(\ga) \mbox{ is an outer corner  of } \gl\}$$
and  
$$X_{-\ga} = \{\gl\in X| \fB(\ga) \mbox{ is an inner corner  of } \gl\}.$$  
For 
$\gl\in  X_{\ga}$,  $t_{ \ga}(\gl) \in X_{-\ga} $ is obtained by adding 
$ \fB(\ga)$ to $\gl$.  Similarly for 
$\gl\in  X_{-\ga}$,  $t_{- \ga}(\gl) \in X_{\ga} $ is obtained by deleting 
$ \fB(\ga)$ from $\gl$.  
This defines mutually inverse bijections  $ X_{\ga} \rl X_{-\ga} .$

\bl \label{rjs} 
If $\ga = \gep_i-\gd_{j}$, the following are equivalent
\bi \itema $ \fB(\ga)$ is an inner corner  of $\gl$, that is $\gl\in X_{-\ga}$.
\itemb $\gl_{n+1-i} = j$ 
and $\gl'_{j} = n+1-i$ .
\ei
If these conditions hold then 
\be \label{qef}  t_{-\ga}(\gl) = (\gl_1,\ldots , \gl_{n-i}, \gl_{n+1-i}-1,\gl_{n+2-i},\ldots, \gl_n).\nn\ee
\el
\bpf  The equivalence of (a) and (b) follows since $ \fB(\ga)$ is an inner corner  of $\gl$ iff $ \fB(\ga)$ is the rightmost entry in row $\gep_i$ and the highest entry in column $\gd_{j}$ of the Young diagram $\gl$.
The formula for $t_{-\ga}(\gl)$ follows from the definitions.\epf 
\noi Similarly we have
\bl \label{rjx} 
If $\ga = \gep_i-\gd_{j}$, the following are equivalent
\bi \itema $ \fB(\ga)$ is an outer corner  of $\gl$, that is $\gl\in X_{\ga}$.
\itemb $\gl_{n+1-i} = j-1$ and $\gl'_{j} = n-i$.
\ei
If these conditions hold then 
\be \label{qhf} t_{\ga}(\gl) = (\gl_1,\ldots , \gl_{n-i}, \gl_{n+1-i}+1,\gl_{n+2-i},\ldots, \gl_n).\ee
\el \noi 
The partition 
$ \boldsymbol \gl  = (m,1^{n-1})$ plays an important role.

 \subsection{Words}\label{1bs}

 \subsubsection{Young diagrams  and words}\label{bsa2}
To a Young diagram $\gl$ in $X$, we associate a word $w=\ttw(\gl)$ of length $m+n$
in the alphabet $\{\ttr, \ttd\}$.  The upper border of $\gl$ in $X$ is a connected sequence of $m+n$ line segments of unit length starting at the 
top left corner and ending at the bottom right corner of
 {\bf R}.  Each line segment carrries the label
$\ur$ or $\ud$, according to whether it moves to the right or down.  Then $w$ is the sequence of labels.  Let $\W$ be the set of words consisting of   
$m$ $\ur$'s and $n$ $ \ud$'s.  There is a bijection  $\ttw:X\lra \W$ given by
$\gl \lra\ttw(\gl)$. We need  a rotation operator $R$ for words. 
 For a word 
\be \label{inb} w= x_1x_2\ldots x_{m+m} \in \W, \mbox{ define } R(w)= x_2\ldots x_{m+m}x_1.\ee

\subsubsection{Diagrams of type $\widetilde{A}(n-1|m-1)$} \label{dta} There is another natural way in which word arise. Let $r=m+n-1$  and $n\neq m$.
A  Coxeter diagram of type $\widetilde{A}_r$ consists of a  
simple closed polygon  with $n+m$ edges, connected by $n+m$  nodes.  Modify the diagram so that it contains  an even, positive number of 
 grey nodes.. 
Draw the diagram with nodes labelled $1,2,\ldots r$, from left to right below the extending node labelled 0.    Moving in a anti-clockwise direction label each edge with a symbol $\ttr$ or $\ttd$ according to the following rules.
\bi \itema The edge following node 0 has label $\ttr$,
\itemb  If node $i$ is white then the labels on the edges before and after node $i$ are the same.
\itemc If node $i$ is grey then the labels on the edges before and after node $i$ different.   
\ei
Then the modified diagram is  a {\it diagram of type} $\widetilde{A}(n-1|m-1)$
if at  the end of this process  the total number of edges labelled $\ttr$ is $m$ or $n$. 
However if this number is $n$, start again assuming in place of (a) that the edge following node 0 has label $\ttd$.  Thus we assume (b), (c) and 
\bi \itemd The total number of edges labelled $\ttr$ is $m$.
\ei 
Now starting at node $i$  and  moving in an anti-clockwise direction we obtain a word $w_i$ in  $\ttr, \ttd$.  For some  Young diagram $\gl_{i}$, we have $w_i=\ttw(\gl_{i})$.
 Thus $R^i(w_0) = w_i$.
\bexa \label{iob} {\rm Let $(n,m)=(2,3)$. The  words $w_0$ associated to the first, second and third diagrams below are $\ttd \ttd\ttr \ttr \ttr$,  $\ttr \ttr \ttr \ttd \ttd$ and $\ttr \ttd \ttr  \ttd \ttr $ respectively.  For the first diagram we have $w_1= \ttd\ttr \ttr \ttr\ttd$, $w_2= \ttr \ttr \ttr  \ttd \ttd $, $w_3 =  \ttr \ttr \ttd \ttd\ttr $ and  $w_4= \ttr\ttd \ttd \ttr \ttr.$   
\vspace{1cm}

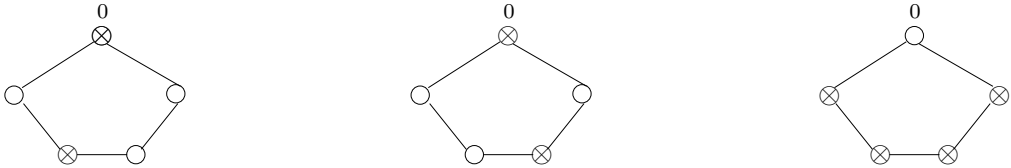
\Bc
\setlength{\unitlength}{0.9cm}
\begin{picture}(10,1)(-0.5,1.70)
\thinlines
  \linethickness{.09mm}
\put(-2.3,0.82){{$\otimes$}}
\put(-1.13,0.93){\circle{0.27}}
\put(-1.8,2.58){$\otimes$}
\put(-0.54,1.83){\circle{0.27}}
\put(-1.8,2.58){$\otimes$}
\put(-2.93,1.81){\circle{0.27}} 
\put(-0.47,1.94){\line(-5.2,3){1.09}}
\put(-2.81,1.92){\line(5,3.2){1.06}}
\put(-2.25,1.01){\line(-4,5){0.54}}
\put(-1.05,1.01){\line(4,5){0.54}}
\put(-1.7,2.95){$\scriptstyle 0$}
\put(-2,.93){\line(1,0){0.74}}
\thinlines
  \linethickness{.09mm}
\put(3.87,0.93){{\circle{0.27}}}
\put(4.7,0.82){{$\otimes$}}
\put(4.2,2.58){$\otimes$}
\put(3.07,1.81){\circle{0.27}} 
\put(5.46,1.83){\circle{0.27}}
\put(4.3,2.95){$\scriptstyle 0$}
\put(5.53,1.94){\line(-5.2,3){1.09}}
\put(3.19,1.92){\line(5,3.2){1.06}}
\put(3.75,1.01){\line(-4,5){0.54}}
\put(4.95,1.01){\line(4,5){0.54}}
\put(4.0,.93){\line(1,0){0.74}}
\thinlines
  \linethickness{.09mm}
\put(10.3,2.95){$\scriptstyle 0$}
\put(9.7,0.82){{$\otimes$}}
\put(10.7,0.82){$\otimes$}
\put(10.37,2.69){\circle{0.27}}
\put(8.95,1.7){{$\otimes$}} 
\put(11.46,1.7){$\otimes$}
\put(11.53,1.94){\line(-5.2,3){1.09}}
\put(9.19,1.92){\line(5,3.2){1.06}}
\put(9.75,1.01){\line(-4,5){0.54}}
\put(10.95,1.01){\line(4,5){0.54}}
\put(10.0,.93){\line(1,0){0.74}}
\end{picture}
\Ec 
\vspace{1cm}
}\eexa \noi
Using words it is easy to describe the equivalence classes for $\sim$ in Subsection 
\ref{tgs} in a way that is very convenient for the proofs of Proposition \ref{crb} and
the main result of \cite{M24}.  
The analog of  odd reflections for 
 words is also easy to describe.  For a positive morphism an occurence of the consecutive letters $\ttd\ttr$ is repaced by $\ttr\ttd$.  More details of this process, which can be written   $p_{\ga}:\W(\ga)\lra \W(-\ga)$ are given 
before Lemma \ref{32s}.

 \subsubsection{Shuffles}\label{shf}
Let 
\be\label{Imn}  I_0= \{1, \ldots, n\}, \quad I_1=\{1' , \ldots , m'  \}\nn\ee and 
let
$Sym_I$ be the group of permutations of $I= I_0 \cup I_1$.    We refer to elements of $I_0 $ (resp. $I_1$) as {\it unprimed} (resp. {\it primed}) integers.
If $\gs \in Sym_I$  the {\it one line notation} for the permutation $\gs$ is
\be \label{olnot}   { \boldsymbol \gs} 
= ( \gs(1), \ldots,\gs(n),\gs(1'),  \ldots \gs(m') ). \ee
The set of 
{\it shuffles} is
\[ \Sh = \{ \sigma \in Sym_I| 1, \ldots, n\;\; {\rm and} \;\; 1',\ldots,  m' \;\; {\rm are \; subsequences \;of}\;\;{ \boldsymbol \gs}  \}. \]
\noi 
From any word $\ttw\in \W$ we obtain a shuffle $\gs=\sh(\ttw)$.  Reading $\ttw$ from left t right, $\gs$ is obtained by replacing each occurrence of $\ttd$ by the smallest available unprimed integer and  each occurrence of $\ttr$ by the smallest available primed integer.  
We have 
\bl \label{xmn} There is a bijection between $\gz:\Sh \lra X$ defined as follows: for a shuffle  $\gs$ we
draw a path in {\bf R}, where  the $k^{th}$ step is $\ud$ if $\gs(k) \le n$
and $\ur$ otherwise. Then the Young diagram $\gl= \gz(\gs)$ consists of the boxes below the path.
\el 
\bpf See  \cite{BN}  Lemma A.1   and Example A.2.\epf \noi 
\bl We have a commutative diagram of bijections.
\be \label{bby2}
\xymatrix@C=2pc@R=1pc{
&\Sh\ar@{<-}^{\sh}[ddr] \ar@{->}_{\gz}[ddl] &
&\\ \\
X \ar@{-}[rr]_\ttw &&
\W\ar@{<-}[ll] &}
\ee
\el \noi 
Mostly we use partitions and words,  but for example from a  shuffle it is very easy to write the names of the simple roots of the corresponding Borel subalgebra, see  \eqref{pgs} and Example \ref{E4}.

 \subsubsection{Borel subalgebras of $\fgl(n|m)$ and odd reflections}\label{0}

Let ${{\stackrel{{\rm o}}{\cB}}}$ be the set of Borel subalgebras $\fbo$ of $\fgo=\fgl(n|m)$ whose even part consists of the upper triangular matrices in $\fgo_0$. 
We can use shuffles to index Borels in ${{\stackrel{{\rm o}}{\cB}}}$  and list their simple roots.  There is a bijection $\Sh\lra  {{\stackrel{{\rm o}}{\cB}}}$, 
$\gs\lra {\stackrel{{\rm o}}{\fb}}{}(\gs)$, where for $\gs\in \Sh$  the corresponding Borel subalgebra 
${\stackrel{{\rm o}}{\fb}}{}(\gs) \in {{\stackrel{{\rm o}}{\cB}}}$ has set of simple roots 
\be \label{pgs} 
{\stackrel{{\rm o}}{\Pi}}(\gs)=\{ \gep_{\gs(i)}- \gep_{\gs(i+1)}|i\in [m+n-1]\}.
\ee  
Thus the distinguished Borel subalgebra ${\stackrel{{\rm o}}{\fb}}{}^\dist$
corresponds to $\gs=1.$  Now  we  have   bijections 
\be \label{bji} \W \rl X \rl \Sh\rl  {\Hbo},\quad \quad w=\ttw(\gl) \rl \gl=\gz(\gs) \rl \gs\rl {\stackrel{{\rm o}}{\fb}}{}(\gs).\ee
If $\ga = \gep_i-\gd_{j}$, we determine the analogs for $\W$ and ${{\stackrel{{\rm o}}{\cB}}}$ of the bijections 
$t_{\pm\ga}: X_{\pm\ga} \lra X_{\mp\ga} $ from Subsection \ref{ypc}. 
Define \be \label{323} 
\W(\ga)= \{w\in \W|\ga \in{\stackrel{{\rm o}}{\Pi}}(\sh(w))\},\quad {\stackrel{{\rm o}}{\cB}}(\ga)=  \{ \fbo(\gl)| \gl \in X, \;\ga\in 
{\stackrel{{\rm o}}{\Pi}}(\sh(\ttw(\gl)))\}.\ee
Let $\fbo\in {\stackrel{{\rm o}}{\cB}}(\ga)$ be a Borel with simple roots ${\stackrel{{\rm o}}{\Pi}}$. 
For $\gb \in \pigo$ we define a root $r_\ga(\gb)$ by
\[ r_\ga(\gb) = \left\{ \begin{array}
  {rcl}
  -\ga  & \mbox{if} & \gb = \ga \\
\ga + \gb & \mbox{if} & \ga + \gb \mbox{ is a root}.\\
\gb & \quad \quad \mbox{otherwise} &
\end{array} \right. \] 
Then $r_\ga(\pigo)$ is the set of simple roots for a Borel subalgebra $r_\ga(\fbo)$. We say $r_\ga(\fbo)$ is obtained from $\fbo$ using the  {\it odd reflection} corresponding to $\ga$. Clearly
$r_{\pm\ga}:{\stackrel{{\rm o}}{\cB}}(\ga)\lra {\stackrel{{\rm o}}{\cB}}(-\ga)$ are mutually inverse bijections.  
If  $\ga =\gep_i-\gd_j$, it follows from \eqref{pgs} and \eqref{323} that 
$w\in \W(\ga)$ (resp. $w\in \W(-\ga)$)  iff $i,j'$ (resp. $j',i$) are consecutive entries in 
$\sh(w)$.   If $i,j'$ appear as entries $k$ and $k+1$ then 
$\ttd, \ttr$ appear in positions $k$ and $k+1$ of $w$, and we obtain $p_\ga(w)\in \W(-\ga)$ by replacing this occurence of $\ttd, \ttr$  by $\ttr,\ttd$.  Thus  $p_\ga$ is a bijection  $p_{\ga}:\W(\ga)\lra \W(-\ga)$ whose inverse we denote by $p_{-\ga}$.

\bl \label{32s} 
The maps from \eqref{bji} restrict to bijections 
\be \label{bj1} X_{\pm\ga} \rl \W(\pm\ga)\rl  {{\stackrel{{\rm o}}{\cB}}}(\pm\ga),\nn\ee
and we have a commutative diagram.
\[
\xymatrix@C=2pc@R=1pc{
X_{\ga}\ar@{<->}[rr] \ar@{->}_{t_\ga} [dd]&&
\W(\ga)\ar@{->}_{p_\ga}[dd]&&
{\stackrel{{\rm o}}{\cB}}(\ga) \ar@{<->} [ll] \ar@{->}^{r_\ga}[dd]&\\ \\
X_{-\ga}\ar@{<->}[rr] &&\W({-\ga})\ar@{}[rr] &&
{\stackrel{{\rm o}}{\cB}}({-\ga})\ar@{<->}[ll] &}
\] 
\el
\bpf The correspondence $\W(\ga)\rl  {{\stackrel{{\rm o}}{\cB}}}(\ga)$ and the commutativity of the right square are  immediate from
\eqref{323}. 
Consider the diagram below. 
On the left (resp.  right) we show the paths determined  $\gl= \gz(\gs)$  
when 
$\ga = \gep_i-\gd_{j}$ is an outer (resp. inner) corner of  $\gl= \gz(\gs)$.  

\Bc
\setlength{\unitlength}{0.8cm}
\begin{picture}(10,3)(-1,-1.50)
\thinlines
\linethickness{0.05mm}
\put(1.39,0.35){$\ga$}
\multiput(1.0,0)(1,0){2}{\line(0,1){1.0}}
 \multiput(1,0)(0,1){2}{\line(1,0){1.0}}
  \linethickness{0.5mm}
\put(8.0,1){\line(1,0){1}}
\put(9.0,1){\line(0,-1){1}}
\put(6.45,0.35){$ \gep_i$}
\put(7.19,0.35){\dots}
\put(8.38,-.75){\vdots}
\put(0.19,0.35){\dots}
\put(1.38,-.75){\vdots}
\put(-.45,0.35){$ \gep_i$}
\put(8.39,-1.35){$ \gd_j$}
\put(1.39,-1.35){$ \gd_j$}
  \linethickness{0.05mm}
\put(8.39,0.35){$\ga$}
\multiput(8.0,0)(1,0){2}{\line(0,1){1.0}}
 \multiput(8,0)(0,1){2}{\line(1,0){1.0}}
  \linethickness{0.5mm}
\put(1.0,0){\line(1,0){1}}
\put(1.0,1){\line(0,-1){1}}
  \end{picture}
\Ec
The remaining statements  follow from Lemma \ref{xmn} and a consideration of these diagrams.
\epf

 \subsubsection{Isomorphisms of groupoids}\label{iog} 
We consider the groupoids from the table.

\[ \begin{tabular}{|c||c||c|} \hline
Groupoid &\;\; Base \;\;& Non-identity morphisms \\ \hline
$\fT(X)$&$\{X_{\pm\ga}|\ga\in\Gd_1^+\}$ &  $\{t_{\pm\ga}|\ga\in\Gd_1^+\}$ 
\\ \hline
$\fT(\W)$&
$\{\W(\pm\ga)|\ga\in\Gd_1^+\}$
 & $\{p_{\pm\ga}|\ga\in\Gd_1^+\}$ 
\\ \hline
$\fT({\stackrel{{\rm o}}{\cB}})$&$\{{\stackrel{{\rm o}}{\cB}}(\pm\ga)|\ga\in\Gd_1^+\}$&$\{r_{\pm\ga}|\ga\in\Gd_1^+\}$ \\ \hline
\end{tabular}\]
It follows from Lemma \ref{32s} that there is a 
$\mathfrak T_{iso}$-equivariant bijection $X\rl {\stackrel{{\rm o}}{\cB}}$.  However more is true.  
There are 
functors 
\bi \itema
$F:\fT_{\iso}\lra \fT(X)$ given by 
$ F({\alpha})=X_\ga$  and 
$F(\gr_{\alpha}) =t_\ga$.
\itemb 
 $G:\fT_{\iso}\lra \fT(\W)$ given by 
$ G({\alpha})=\W(\ga)$  and 
$G(\gr_{\alpha}) =p_\ga$.
\itemc 
 $B:\fT_{\iso}\lra \fT({\stackrel{{\rm o}}{\cB}})$ given by 
$ B({\alpha})={\stackrel{{\rm o}}{\cB}}(\ga)$  and 
$B(\gr_{\alpha}) =r_\ga$. 
\ei Since all we have done is rename the objects and morphisms, it is clear that  $F, G$ and $B$ are isomorphisms of groupoids and we have a commutative diagram of isomorphisms.
\be \label{bby1}
\xymatrix@C=2pc@R=1pc{
\mathfrak  T_{iso}\ar@{<->}_{F} [dd]&&
\mathfrak  T_{iso} \ar@{<->}_{=}[ll]\ar@{<->}_{{G}}[dd]&&
\mathfrak  T_{iso}\ar@{<->}_{B} [dd]\ar@{<->}_{=}[ll]&&&\\ \\
\fT(X)&&\fT(\W) \ar@{<->}^{{\cong}}[ll]&
& \ar@{<->}^{{\cong}}[ll]
\fT({\stackrel{{\rm o}}{\cB}})\ar@{<-}[ll] &}
\nn\ee
\bc \label{APE} We have $\sh(p_\gb(w)) =r_\gb(\sh(w))$ for $w\in \W(\gb)$. \ec \noi
We define $\W(\rda)$, ${\stackrel{{\rm o}}{\cB}}(\rda)$ to be the subsets of $\W, {\stackrel{{\rm o}}{\cB}}$ corresponding to $X(\rda) $. Similarly define  $\W(\rpa)$, $\W(\cda)$,  $\W(\cpa)$, ${\stackrel{{\rm o}}{\cB}}(\rpa)$,  ${\stackrel{{\rm o}}{\cB}}(\cda)$ and  ${\stackrel{{\rm o}}{\cB}}(\cpa)$. Then we use the same superscripts $\pr,  \mr,\opc, \omc$ for bijections between subsets of  $\W$, ${\stackrel{{\rm o}}{\cB}}$ that we used for subsets of $X$.  
\\ \\
To handle several cases simultaneaously,  suppose $H$ is a functor $H:\mathfrak T_{iso}\lra S(\X)$ 
where 

\be \label{BPG}( \X, H)  = (X, F), (\W, G) \mbox{ or  } ({\stackrel{{\rm o}}{\cB}}, B)\ee  as above.  
\subsubsection{ The equivalence relation} \label{tgs} 
In the Introduction we have already introduced an equivalence relation $\sim$  on
$X\ti \Z$.   If $\X$ is as in \eqref{BPG} (and $\X\neq X$), then using the bijections from Subsection \ref{iog},  we obtain a compatible equivalence relation on 
$\X\ti \Z$ which we also  denote by $\sim$.

\bl \label{rogs} \bi \itemo 
\itema $\gl \in X(\rda)$ iff $\ttw(\gl)= \ttx\ttd$  $($for some sub-word $\ttx)$.   In this case $\ttw(\gl^\mr)= \ttd
\ttx.$ 
\itemb  $\gl \in X(\rpa)$ iff $\ttw(\gl)= \ttd\ttx.$ In this case $\ttw(\gl^\pr)= \ttx\ttd
.$ 
\itemc
$\gl \in X(\cda)$ iff 
$\ttw(\gl)= \ttr\ttx.$ In this case $\ttw(\gl^\omc)= \ttx\ttr
.$ 
\itemd $\gl \in X(\cpa)$ iff $\ttw(\gl)= \ttx\ttr.$ In this case $\ttw(\gl^\opc)= \ttr\ttx.$ 
\ei
\el 
\bpf Clearly $\gl \in X(\rda)$ iff the last letter in the path (or word) corresponding to $\gl$ is $\ttd$, that is  $\ttw(\gl)= \ttx\ttd$.  Moving $\ttd$ from the end of $\ttw(\gl)$ to the start,   
we obtain a new path with first letter $\ttd$ and with $i+1^{st}$ letter equal to the $i^{th}$ letter of $\ttw(\gl) $.  
The  new path corresponds to $\gl^\mr$.
This shows (a) and the other parts are similar. 
\epf \noi 
Since any word starts (or ends) with the letter $\ttd$ or $\ttr$ we have.
\bc We have disjoint unions \bi\itema $X =X(\rpa) \cup X(\cda).$
\itemb $X =X(\cpa) \cup X(\rda).$
\ei
\ec
\noi
Using words it is easy to describe the equivalence class of $(w_0, k_0)$  where $w_0 = \ttw(\gl)$.  To do this we define inductively words $w_i$ ($i>0$) using the following rules 
\bi
\itema If  $w_i=\ttd x$, then $w_{i+1} = x\ttd $ and $k_{i+1} = k_i -m$.
\itemb If  $w_i=\ttr x$, then $w_{i+1} = x\ttr $ and $k_{i+1} = k_i +n$.
\ei 
Define $\gl_i \in X$ by   $w_{i} =\ttw (\gl_i)$.
 \bl \label{iod}
 \begin{itemize}\itemo
\itema If the first letter in $w_i$ is $\ttr$, then $\gl_{i+1}=\gl_{i}^\omc$.
\itemb If the first letter in $w_i$ is $\ttd$, then $\gl_{i+1}=\gl_{i}^\pr$.
\ei
\noi Thus by Lemma \ref{rogs} $w_{i+1}$ is obtained from 
$w_{i}$ by moving the first letter to the end. 
 \el
\noi
Now each word contains $m$ symbols $\ttr$ and $n$ symbols $\ttd$.  Therefore $w_{m+n} = w_0.$
To find a closed formula for $k_i$ set 
\[ \begin{array}
  {rcl}
 c_j =1, d_j =0 & \mbox{if} & \mbox{the } j^{th}  \mbox{ letter of } \ttw \mbox{ is } \ttr\\
 c_j =0, d_j =1 & \mbox{if} & \mbox{the } j^{th}  \mbox{ letter of } \ttw \mbox{ is } \ttd\end{array},\] then define
\be \label{ips}  k_{i} = k_0 +n\sum_{j=1}^{i}c_j -m\sum_{j=1}^{i}d_j. 
\ee

\bc \label{vgs} \itemo
\bi \itema 
The set $E:=\{(w_{i},k_{i})| 0 \le  i < m+n\}$, is the  equivalence class of $(w_0,k_0)$ in $[\W\ti\Z]$. 
\itemb The equivalence class of 
$(\fbo(w_0),k_0)$ in 
$[\stackrel{{\rm o}}{\cB} \ti \Z]$ is   $\{(\fbo(w_i),k_i))| 0 \le  i < m+n\}.$
\itemc  If $0\le i \le j\le m+n$ and $k_i\equiv k_j \mod mn$ then $i=j$.
\itemd If $(w_{i},k_{i})\in E$ and $m$ divides $k_i - k_0$, then $\sum_{j=1}^{i}c_j$ equals $ 0 $ or $m$. 
\iteme If $w_0= \ttx\ttr$  for some word $\ttx$, $(w_{i},k_{i})\in E$ and $m$ divides $k_i - k_0$, then $\sum_{j=1}^{i}c_j=0$. 
\ei\ec
\bpf If the first letter in $w_i$ is $\ttr$, then $\gl_{i+1}=\gl_{i}^\omc$ and 
$(w_{i},k_{i}) \sim (w_{i+1},k_{i}+n)$ by \eqref{2kr} and this agrees with \eqref{ips}.
Thus $E$ is contained the equivalence class of $(w_0,k_0)$.  However the set $E$ is closed under the basic instances (\eqref{1kr} and \eqref{2kr}) of $\sim$, so the equivalence class can be no larger.  This shows (a) and the proof of (b)  is similar.   Note that  (c)  and (d) follow from \eqref{ips}. If $w_0= \ttx\ttr$, then $ c_{m+n} =1$  and  $k_0= k_{m+n} =  k_{m+n-1} +n $, so $m$ cannot divide $k_i - k_0$, unless $\sum_{j=1}^{i}c_j=0$.   
\epf\noi 
We remark that words are very useful for computations as in the following example. 

\bexa \label{glk}{\rm If  ${\gl= \boldsymbol \gl},$ then  $w=\ttw(\gl)=\ttr^{}  \ttd^{n-1}\ttr^{m-1}\ttd^{}$. Deleting a row we have $w^\mr=\ttw\gl^\mr)=\ttd^{}\ttr^{}  \ttd^{n-1}\ttr^{m-1}$ by Lemma \ref{rogs} (a).  A positive odd reflection gives 
$\mu=t_\ga(\gl^\mr)$ with $\ttw(\mu)=\ttr^{}  \ttd^{n}\ttr^{m-1}.$ Here $\ga=\gep_1- \gd_1$ Then deleting a column we have $\mu^\omc$ (which corresponds to the empty partition) by Lemma \ref{rogs} (c), and we have $\ttw(\mu^\omc)=\ttd^{n}\ttr^{m}$. Now by \eqref{1kr} 
$(\gl,0)\sim (\gl^\mr, m)$. Then 
$(t_\ga(\gl^\mr),m)=(\mu, m)$ 
and  by \eqref{2kr} 
$(\mu,m)\sim (\mu^\omc, m+n)= (\emptyset, m+n)$.  
The same result is obtained by deleting a column first: we have
$\ttw(\gl^\omc)=  \ttd^{n-1}\ttr^{m-1}\ttd\ttr$ by Lemma \ref{rogs}.  Then an  odd reflection using the root $\gb=\gep_n- \gd_m$ yields  
$\phi=r_\gb(\gl^\omc)$, with $\ttw(\phi)= \ttd^{n-1}\ttr^{m}\ttd^{} .$ Then deleting a row we have $\ttw(\phi^\mr)=\ttd^{n}\ttr^{m}$.
 More generally if we carry out the same sequence of operations starting with a partition  $\xi$ containing  ${ \boldsymbol \gl}, $ we obtain the skew diagram $\xi/{ \boldsymbol \gl}$. This fits into the lower left corner of $X$. 
}\eexa

\section{Extending the action of $\mathfrak  T_{iso}$} \label{APG}

\subsection{Actions of groupoids on equivalence classes }\label{itr1} Suppose $\mathfrak T_{iso}$ acts on $Y$ via a functor $H:\mathfrak T_{iso}\lra S(Y)$ and that $\sim$ is an equivalence relation on $Y$.  For $\ga\in\Gd_{iso}$, set $Y(\ga) = H(\ga)$.  Let 
$[Y(\ga)] $ denote the set of equivalence classes of elements in $Y(\ga).$
We do not assume that $Y(\ga)$ is a union of equivalence classes.  Note that $[y]\in
[Y(\ga)] $ iff there exists some $z \in [y] \cap Y(\ga).$  It is not always possible to take 
$z=y$, see for example Equation \eqref{sb1}.

\bl \label{sogs} Suppose that  $y\sim z$ with $y,z \in Y(\ga)$ implies $H(\gr_\ga)(y)\sim H(\gr_\ga) (z) \in Y(-\ga)$.  Then there is a well-defined action of $\mathfrak  T_{iso}$ 
on  the set of equivalence classes $[Y]$ via a functor  $\H:\mathfrak T_{iso}\lra S([Y])$.  
It is given by
$\H(\ga)= [Y(\ga)] $ and 
$$\H(\gr_\ga)([y])  = [H(\gr_\ga)(z)]$$
where $z \in [y] \cap Y(\ga).$
\el


\subsection{A general construction} \label{AP1}  
Assume  $\mathfrak  T_{iso}$  acts on $\X$ as in \eqref{BPG}  and  that $n,m$ are coprime. Our goal is to construct a functor $\H:\mathfrak  T_{iso}\lra S(\cX)$, where $\cX$ is the set of equivalence classes in $\X\ti \Z$ under the equivalence relation from Subsection \ref{tgs}.  First, 
we define an action of  $\mathfrak  T_{iso}$ on $\X\ti \Z$.
This action preserves the  rotation number, so it is enough to describe the restriction to $(\X,k)$. 
Set 
\be \label{sy1}
H_{in+jm } (\ga)= (\X_{ \eta^{i}\nu^{j}\ga},in+jm )
.\ee

\be \label{szx}
H(\gr_\ga)(x, in+jm) = (H(\gr_{\eta^{i}\nu^{j}\ga}
)(x), in+jm) 
\mbox{ for  } x \in \X_{ \eta^{i}\nu^{j}\ga}
.\ee
This is well-defined since $n,m$ are coprime and 
$\eta, \nu$ have order $m, n$ respectively. Next set
\be \label{vzx}
H(\ga) = \bcu_{k\in \Z}  
H_{k} (\ga) \subset \X\ti \Z\ee and let $\cX_\ga$ be the set of equivalence classes in $H(\ga) $. 
Equations \eqref{szx}
and  \eqref{vzx}
define a functor 
$\mathfrak  T_{iso}\lra \cS(\X\ti \Z).$

\bh \label{ky} 
  Assume that   the following conditions holds.

\bi \itema
If $\gb\in \Gd_{1}^+$, $\gl\in \X_{\gb}\cap  \X(\rda)$ and 
${\gl}^\mr \in \X_{\nu\gb}$, then ${H(\gr_{\gb})(\gl)}\in  \X(\rda)$ and   
\be \label{kz} ({H(\gr_{\gb})(\gl)})^\mr = H(\gr_{\nu\gb})({\gl}^\mr).\ee 
\itemb  
If $\gb\in \Gd_{1}^+$, $\gl\in \X_{\gb}\cap  \X(\cda)$ and 
${\gl}^\omc \in \X_{\eta\gb}$, then ${H(\gr_{\gb})(\gl)}\in  \X(\cda)$ and   
\be \label{kz1} ({H(\gr_{\gb})(\gl)})^\omc = H(\gr_{\eta\gb})({\gl}^\omc).\ee 

\ei\eh
\bl Assume Hypothesis \ref{ky} holds. \bi \itema
If  $\gb\in \Gd_{1}^+$,  $\mu\in  \X_{-\nu\gb}\;\cap \in \X(\rpa)$ and  ${\mu}^\pr \in \X_{-\gb}$, then  ${H(\gr_{-\nu\gb})(\mu)}\in  \X(\rpa)$ and  
\be \label{kf} ({H(\gr_{-\nu\gb})(\mu)})^\pr = H(\gr_{-\gb})({\mu}^\pr).\ee
\itemb
If  $\gb\in \Gd_{1}^+$,  $\mu\in  \X_{-\eta\gb}\;\cap \in \X(\cpa)$ and  ${\mu}^\opc \in \X_{-\gb}$, then  ${H(\gr_{-\eta\gb})(\mu)}\in  \X(\cpa)$ and  
\be \label{kf1} ({H(\gr_{-\eta\gb})(\mu)})^\opc = H(\gr_{-\gb})({\mu}^\opc).\ee
\ei
\el 
\bpf For (a), set $\mu = ({H(\gr_{\gb})(\gl)})^\mr = H(\gr_{\nu\gb})({\gl}^\mr)$ and solve  
 for $\gl$ in two different ways.
\epf
\noi 
Given this data,  we explain how to define an action of $\mathfrak T_{iso}$ on the set of equivalence classes $\cX=[\X\ti \Z]$ under the equivalence relation generated by \eqref{1kr} 
and \eqref{2kr}. 
\bt \label{shg}  If the functor $H$ satisfies Hypothesis \ref{ky}, then 
\bi \itema There is a functor $\H:\mathfrak  T_{iso}\lra \cS(\cX)$ given by $\H(\ga) = \cX_\ga $ 
and
\be \label{syx}
\H(\gr_\ga)[\gl, in+jm] =
 [H(\gr_{\eta^{i}\nu^{j}\ga})(\gl),in+jm]
\mbox{ for  } (\gl, in+jm)\in (X_{ \eta^{i}\nu^{j}\ga},in+jm )
.\ee

\itemb 
The injective map 
$\X\lra \cX,$ $\gl\lra [\gl,0]$ is compatible with the 
 $\mathfrak T_{iso}$-action.  If 
$\gl \in X_{\ga}$, then $\H(\gr_{\ga})[\gl,0]= [H(\gr_{\ga})(\gl),0].$
\ei
Thus the action of 
 $\mathfrak T_{iso}$ on $\X$ extends to an action on $[\X\ti Z].$
\et
\bpf (a) Since $(\gl,in+ jm)\sim  ({\gl}^\mr,in+(j+1)m)$ and $(\gl,in+ jm)\sim  ({\gl}^\omc,(i+1)n+m)$ , to  show the action in \eqref{syx} is well-defined, by Lemma \ref{sogs} we need two conditions.

\bi 
\itemi If 
$\gl \in X_{ \eta^{i}\nu^{j}\ga}$ and  ${\gl}^\omc\in X_{  \eta^{i+1}\nu^{j}\ga}$ then
$$(H(\gr_{ \eta^{i} \nu^{j}\ga})(\gl),in+jm)\sim (H(\gr_{   \eta^{i+1}\nu^{j}\ga})({\gl}^\omc),(i+1)n+jm).$$
\itemii If 
$\gl \in X_{ \eta^{i}\nu^{j}\ga}$ and  ${\gl}^\mr\in X_{  \eta^{i}\nu^{j+1}\ga}$ then
$$(H(\gr_{ \eta^{i} \nu^{j}\ga})(\gl),in+jm)\sim (H(\gr_{  \eta^{i}\nu^{j+1}\ga})({\gl}^\mr),in+(j+1)m).$$

\ei For (i), set  $k= in+jm$ and  $\gb= \eta^{i}{ \nu^{j}\ga}$, then
using the definition of $\sim$ 
 and \eqref{kz1}
\by \label{a12}   
(H(\gr_{\gb})(\gl),k)
&\sim& 
({H(\gr_{\gb})(\gl)}^\omc,k+n)
\nn\\&=& 
(H(\gr_{\eta\gb})({\gl}^\omc),k+n).
\nn
\ey  as required.  The proof of (ii) is similar using \eqref{kz}. 
In statement (b), injectivity of the map  $\X\lra [\X\ti Z],$  follows from Corollary  \ref{vgs} (c) and the other statement follows from \eqref{syx} with $i=j=0$.\epf

\subsection{Application to Young diagrams, words and Borel subalgebras} \label{AP2} 
 {\bf Notation: }We assume $w, \gl$ and $\gs$ are related as in \eqref{bji}.
\bl \label{rhy} If
$\gl\in  X(\rda)\cap  X_{\gb}$, then ${\gl}^{-_\ttr} \in 
X_{\nu\gb}$,  ${t_{\gb}(\gl)}\in X(\rda)$ and  
\be \label{kzt} {t_{\gb}(\gl)}^\mr = t_{\nu\gb}({\gl}^\mr).\ee
\el
\bpf Suppose $\gb= \gep_i-\gd_{j}\in \Gd_{1}^+$ and set $t_\gb=   F(\gr_{\gb})$. 
Since $\gl_1 = m$, $\gl$ has no outer corners in row $\gep_{n}$. Thus 
 $i<n$.
The equivalence of (a) and (b) in  Lemma \ref{rjx} implies 
$\nu\gb= \gep_{i+1}-\gd_j$  is an outer corner of ${\gl}^\mr$ and the rest follows from \eqref{qhf}.  
\epf 
\noi The analog of \eqref{kzt} for words is 
 \be \label{kut} {p_{\gb}(w)}^\mr = p_{\nu\gb}({w}^\mr).\ee
To see this directly let $w=u\ttd \ttr v \ttd$ for sub-words $u,v$. Then 
$w^\mr = \ttd u \ttd \ttr v$, 
 $p_{\gb}(w)=u\ttr \ttd v \ttd$ and 
$$p_{\gb}(w)^\mr = p_{\nu\gb}(w^\mr)=\ttd u\ttr  \ttd v.$$
\noi Recall the bijection 
$\gz:\Sh \lra X$ given in Lemma \ref{xmn}. 
 Suppose  $\gl=\gz(\gs) \in X(\rda)$ 
We determine the shuffle that corresponds to  
${\gl}^\mr$. 
 Let  ${ \gs}^\mr$ be the   permutation whose one-line notation is
$${ \boldsymbol \gs}^\mr=(1,  \nu^{-1}\gs(1), \ldots,\nu^{-1}\gs(n),\nu^{-1}\gs(1'),  \ldots \nu^{-1}\gs(m+n-1) ).$$

\bl \label{ymn} We have 
\bi \itema
${\gl}^\mr =\gz({\sigma}^\mr).$
\itemb $\sh(w)^\mr= \sh(w^\mr)$, 
\itemc  Suppose $\gb =\gep_i-\gd_j$ and  
$r_{\gb}(\gs) = (i,j')\gs$ is obtained from $\gs$ using the { odd reflection} $r_{\gb}$.
Then  
\be \label{e6z}({r_\gb\gs})^\mr= r_{\nu\gb}({\gs}^\mr).\ee 
\ei\el
\bpf   Since the first entry in ${ \boldsymbol \gs}^\mr$ is 1, the first step in the path corresponding
to $\gz({\sigma}^\mr)$ is $\ud$.  
For $i\ge 0$, step $i+1$ of this path is the same as 
step $i$ of the path for $\gl$.  
This proves (a).
If $\ttw(\gl)= \ttx\ttd$ and $\ttw(\gl^\mr)= \ttd
\ttx$ as in Lemma \ref{rogs}, the same reasoning 
 gives (b). 
  For (c),  we use  (b), \eqref{kut} and Corollary \ref{APE}.  Thus 
$$r_{\gb}(\sh(w))^\mr = \sh (p_{\gb}(w))^\mr = \sh(p_{\gb}(w)^\mr) = \sh(p_{\nu\gb} (w^\mr)
=r_{\nu\gb}(\sh(w^\mr)).$$ 
\epf \noi 
We mention the analogs of Equation \eqref{kzt}, \eqref{kut} and \eqref{e6z} that correspond to deleting a column from a Young diagram.
\be \label{kt} {t_{ \gb}(\gl)}^\omc = t_{\eta\gb}({\gl}^\omc),\quad {p_{\gb}(w)}^\omc = p_{\eta\gb}({w}^\omc),\quad ({r_{\gb}}\gs)^\omc= r_{\eta\gb}({\gs}^\omc).\ee

\noi 
\bl \label{rhk} The functor $H$ satisfies Hypothesis \ref{ky}.
\el
\bpf Part (a) of the Hypothesis follow from Equations \eqref{kzt}, \eqref{kut}, \eqref{e6z} and a short computation.  In  the Young diagram case under the hypotheses of Lemma \ref{rhy}, we have 
$(\gl, k)\sim (\gl^\mr, k+m)$.  Suppose that $k =in+jm$ and set and  $\xi= \eta^{i}{ \nu^{j}\gb}$.  
Then by \eqref{kzt}
\be \label{kzx}F(\gr_\xi)  ({\gl}^\mr,k+m) = (t_{\nu\gb}({\gl}^\mr),k+m )=
({t_{\gb}(\gl)}^\mr,k+m) \sim   ({t_{\gb}(\gl)},k)=
F(\gr_\xi)(\gl,k)
\ee
 Part (b) follows from   \eqref{kt}.  
\epf

\bc \label{sz} \bi \itemo
\itema If $(\X, H)$ are as in 
\eqref{BPG}, 
there is a well-defined functor $\H:\fT_{\iso}\lra \fS[X\ti \Z]$.  

\itemb 
\noi 
We have a commutative diagram.
\be \label{bby7}
\xymatrix@C=2pc@R=1pc{
&&\mathfrak  T_{iso}\ar@{<->}_{\F} [dd]&&
\mathfrak  T_{iso} \ar@{<->}_{=}[ll]\ar@{<->}_{{\G}}[dd]&&
\mathfrak  T_{iso}\ar@{<->}_{\B} [dd]\ar@{<->}_{=}[ll]&&&\\ \\
&&\cS[X \ti \Z ]&&\cS[\W \ti \Z ] \ar@{<->}^{{\cong}}[ll]&
& \ar@{<->}^{{\cong}}[ll]
\cS[{\stackrel{{\rm o}}{\cB}}\ti \Z ]\ar@{<-}[ll] &}
\nn\ee
\ei
\ec

\bpf  This follows directly from Lemma \ref{rhk} and Theorem  \ref{shg}.
\epf \noi  We study the action of $\fT_{\iso}$ on $[\Hbo\ti \Z]$ (resp. $[X\ti \Z]$) 
further in Section \ref{akmls} (resp.    \cite{M24}). 
\br{\rm  
Injectivity fails in Theorem \ref{shg} (b) if $(n,m)$ are not 
coprime.  Suppose $M =ma =nb$, where $a<n$ and $b <m$ and let $\gl= m^n$.  Deleting $a$ rows from $\gl$ we obtain $\mu = m^{n-a}$ and deleting $b$ columns we obtain $\xi = (m-b)^{n}$.  By \eqref{1kr} and \eqref{2kr} we have
$(\mu,0)\sim ({\gl} ,-M) \sim ({\xi} ,0)$ but
$\mu \neq  {\xi}.$ 
}
\er
 \subsection{$[X\ti \Z]$ as a graded poset}\label{gax}
We can regard $[X\ti \Z]$ as a  poset where $[{\mu},\ell] \le [{\gl},k]$ means that $[{\gl},k]$ can be obtained from $[{\mu},\ell]$ by applying a sequence of positive morphisms.  
We give the pair $({\gl},k)$ a degree by setting
\be \label{e67n}\deg({\gl},k)=|\gl| + k\in \Z.\ee  This degree passes to equivalence classes. For example

\bl \label{rjg1} 
$\deg({\gl}^\mr,k+m)=\deg(\gl,k)$ if $\gl_1=m$. 
\el  \bpf This is clear since $\gl$ contains $m$ more boxes than ${\gl}^\mr$. 
\epf  \noi Thus 
$[X\ti \Z]= \bcu_{d\in \Z} [X\ti \Z]_{d}$ is a graded poset where 
$[X\ti \Z]_d = \{[{\gl},k]|\deg[{\gl},k] =d\}.$  Set $(X\ti \Z)_d = \{({\gl},k)|\deg({\gl},k) =d\}.$  
  For any $\gl\in X$ and $d\in \Z$ there is a unique integer $k$ such that $\deg({\gl},k)=d$, by \eqref{e67n}.  Thus $|(X\ti \Z)_d | = |X|$.  The set $(X\ti \Z)_d $ is a union of  equivalence classes and each equivalence class consists of $m+n$ pairs by Corollary \ref{vgs}.  Hence 
for any $d\in \Z$, we have 
$$|[X\ti\Z]_d|=
\frac{1}{m+n} \left(\begin{array}{c} m+n\\ n
\end{array} \right). $$ 
\noi 
It is clear that $[X\ti \Z]$ enjoys a periodicity property: there is a bijection $[X\ti \Z]_{d}\rl [X\ti \Z]_{d+mn}$,  $[{\gl},k]\rl [{\gl},k+mn]$. 
 \subsection{Pseudo-corners}\label{gbx}
For simplicity we assume $(\X,H)=(X,F) $ in this and the following Subsection, though everything we say has an analog for the other cases of \eqref{BPG}.
 The set of {\it reduced diagrams} is 
$$X^{\str}= \{\gl \subseteq {\bf R} |  \gl_n = \gl'_m =0\}.$$ 
If ${ \boldsymbol \gl} \subseteq \gl$,
we call $\gep_n-\gd_1$   a   {\it outer pseudo-corner} of  $\gl$, and if 
$\gl \in X^{\str}$ we call $ \gep_1-\gd_m  $   a {\it  inner  pseudo-corner} of  $\gl$.  These are not outer  or inner corners, but behave in some ways as if they were.  Pseudo-corners account for the cover relations in the poset $[X\ti \Z]$ that do not arise from $X$  and Theorem \ref{shg} (b).  
It is easy to see the following.
\bl \label{rey} \bi \itemo
\itema  If 
$\ga =\gep_n-\gd_1$ 
and ${ \boldsymbol \gl} \subseteq  \gl$,
 then $\ga$ is not an outer corner of $\gl$, but 
$\nu\ga = \gep_1-\gd_1$ is an outer corner of 
${\gl}^\mr$ and $\eta^{}\ga=\gep_n 
-\gd_m$ is  an outer corner of ${\gl}^\omc.$
\itemb  
If $\ga =\gep_1-\gd_m$ and $\gl \in X^{\str}$, then $\ga$ is not an inner corner of $\gl$, but $\nu^{-1}\ga=\gep_n 
-\gd_m$ is  an inner corner of ${\gl}^\pr$ and $\eta^{-1}\ga=\gep_1 
-\gd_1$ is  an inner corner of ${\gl}^\opc.$
\ei
\el  \noi
We pay special attention to pseudo-corners in Corollary \ref{sz}.  Suppose (a) holds in Lemma \ref{rey} and set 
$\gb=\nu\ga = \gep_1-\gd_1$    and $\gc =\eta^{}\ga=\gep_n 
-\gd_m$. 
We have $(\gl,k)\sim({\gl}^\mr,k+m) \sim({\gl}^\omc,k+n).$
 Hence 
\by \label{1no} \F(\gr_\ga)[\gl,k]&=&\F(\gr_\ga)[{\gl}^\mr,k+m]= [t_\gb{\gl}^\mr,k+m]\\&=&\F(\gr_\ga)[{\gl}^\omc,k+n]= [t_\gc{\gl}^\omc,k+n]\nn
.\ey
Now suppose (b) holds in Lemma \ref{rey} and set 
$\gb=\nu^{-1}\ga=\gep_n 
-\gd_m$ and  $\gc=\eta^{-1}\ga=\gep_1 
-\gd_1$. We have 
$(\gl,k)\sim({\gl}^\pr,k-m) \sim({\gl}^\opc,k-n)$ and 
\by \label{2no}\F(\gr_{-\ga})[\gl,k]&=&\F(\gr_{-\ga})[{\gl}^\pr,k-m]= [t_{-\gb}{\gl}^\pr,k-m].
\\&=&\F(\gr_{-\ga})[{\gl}^\omc,k-n]= [t_{-\gc}{\gl}^\omc.k-n]\nn\ey

 \subsection{Hasse diagrams}\label{gux}
We consider the Hasse diagram 
 for the periodic $\Z$-graded poset $[X\ti \Z]$.  The Hasse diagram shows only the arrows of the Cayley graph $\Gc[X\ti \Z]$ corresponding to positive morphisms.  In poset language there is  an arrow $e\lra e'$  precisely when  $e$ is covered by $e'$.
Before giving an example we prove a result, Theorem \ref{vss} which greatly simplifies our task.  
We can define a new equivalence relation $\ap$ on $X\ti \Z$ 
by requiring that only one of Equations \eqref{1kr} or \eqref{2kr} holds. 
To be defnite we adapt \eqref{1kr}.  Thus let  $\ap$ is the equivalence relation generated by 
$(\gl,k)\ap ({\gl}^\mr,k+m)$ if  $\gl\in X(\rda)$. The proof of Corollary \ref{sz} shows that 
$\fT_{\iso}$
acts on the set of equivalence classes 
$[X\ti \Z]_\ap.$ The connected component of $\Gc[X\ti \Z]_\ap.$  containing $[\emptyset, 0]$ is 
$\Gc[X\ti m\Z]_\ap.$ 
\\ \\
It will be convenient to work with words rather than partitions.  
Recall that in Corollary \ref{vgs}, $E$ denotes the equivalence class of  $(w_0,k_0)$ under $\sim.$ 
We determine the decompoosition of $E$ into $\ap$-equivalence classes.  This depends on maximal sequences of consecutive entries in a word all equal to $\ttd$.  
By replacing $w$ by a cyclic rotation, we may  assume that $w = w_0$ has the form $w= \ttx\ttr$ for some word  $\ttx.$  This ensures that $w$  has no "wrap-around " sequences of the above form.  
We need some index sets: for $j  \in  \Z_m = \{0,1,\ldots m-1\},$ set  
\be  \label{jps}  \cE_j = \{i \in \{0,1,\ldots m+n-1\} | k_i -k_0 \equiv nj \mod m\}.\ee
\bt \label{kps}  If $E_{j} = \{(w_i,k_i)| i \in \cE_{j}\},$ then  $E = \bcu_{i=0}^{m-1} E_j$ is the decomposition of $E$ into $\ap$-equivalence classes.
\et \noi \bpf Define a  map $f:\{0,1,\ldots m+n-1\} \lra \Z_m$ so that $f(i)$ is minimal with $i\in \cE_{f(i)}$. 
By \eqref{ips}  
\be \label{ops}  k_{i+1} = k_i +nc_{i+1} -md_{i+1}. 
\ee
Note that the $i+1^{st}$ entry in $w_0$ equals the first entry in $w_i$.  There are two cases
\bi \itema  $w_i = \ttr\ttx$ for some word $\ttx$ and $c_{i+1} =1.$  Then 
 $k_{i+1} = k_i +n$ and so $f(i+1) =f(i) +1$.
\itemb $w_i = \ttd\ttx$ for some word $\ttx$ and $c_{i+1} =0$.  Then 
 $k_{i+1} = k_i \mod m$ and thus $f(i+1) =f(i).$
\ei  Also $f$ is onto since $w$ contains $m$ letters equal to $\ttr.$ Hence each $E_j$ is  non-empty.  
It follows from (b) that each $E_{j}$ is  a union  of $\ap$-equivalence classes.  Also if $ (w_i,k_i)\ap
 (w_\ell,k_\ell)$ then $m$ divides $k_\ell-k_i$.  Thus the $E_j$ are precisely  the  $\ap$-equivalence classes.
\epf
\bl \label{vgt} 
\bi \itemo  \itema
$E$  contains a representative of the form $(w,\ell)$
such that $m$ divides $\ell$.
\itemb The set of representatives as in (a) forms an  equivalence class of   under $\ap.$ 
\itemc  Suppose that $(w_0, k_0)$ has the property in (a), and  
$[w_0, k_0] \lra p [w_0, k_0]$ is a positive edge  of the Cayley graph $\Gc[\W\ti \Z]$  then $p$ has a representative of the form $(w_0, k_0)\lra (p_\gc (w_0), k_0)$ or $(w_0, k_0)\lra (p_\gc (w_0), k_0+m)$. 
\ei \el
\bpf  Given $(w_0,k_0)$ solve $k_0+nj \equiv 0 \mod m$ for $j$. Then if $i \in \cE_j$, $k_i\equiv 0$ mod $m$.   This gives (a) and (b) follows from \eqref{ips}.  
\\ \\
The  hypotheses in (c) implies that for some $j$ and  root $\gb\in\Gd_1^+$, $w_j=R^{j}(w_0) \in W(\gb)$, see \eqref{inb} for notation. Then 
$\gb$ corresponds to a consecutive pair $\ttd \ttr$ in $w_j$ in entries $k, k+1$: 
$$ w_j=y_1y_2, \ldots y_{m+n}$$
 with $y_{k} = \ttd,     y_{k+1}=\ttr$.  In $w_0=R^{-j}(w_j),$ there is another consecutive pair $\ttd \ttr$  in entries $k+j, k+j+1$. Unless 
$w=\ttr u \ttd $ 
and $k+j\equiv m+n \mod m+n$,  this pair yields a root $\ga\in\Gd_1^+$ such that 
$w_0 \in W(\ga)$.  Let $w'_0=p_\ga(w_0)$ and $w'_j=p_\gb(w_j)$.  Then $w'_j=R^{j}(w'_0)$  and $(w'_0, k_0)\sim (w'_j,k_j)$.   Then $p$ is represented by $(w_0, k_0) \lra  ({p_{\ga}(w_0)},k)$.  
 In the exceptional case where 
$w=\ttr u \ttd $ 
and $k+j\equiv m+n \mod m+n$,    ${ \boldsymbol \gl} \subseteq  \gl$  and 
$\ga =\gep_n-\gd_1$ is an outer pseudocorner of $\gl$ as in   Lemma \ref{rey}.  Then  by Equation \eqref{1no}, $p$ is represented by $(w_0, k_0) \lra  ({p_{\gc}(w_0^\mr)},k+m)$ where $\gc = \nu \ga$. 
\epf 
\bt \label{vss} 
There is  a  
$\fT_{\iso}$-equivariant bijection $[X\ti m\Z]_\ap \lra [X\ti \Z]$ given by
$$[\gl,  m\ell]_\ap \lra [\gl,  m\ell].$$
\et
\bpf The map is well-defined because $e\ap e'$ implies $e\sim e'$ and is clearly 
$\fT_{\iso}$-equivariant.  It is onto by Lemma \ref{vgt}. To show it is injective, suppose $ m\ell =k_0$ and the pairs $(\gl,  m\ell)$ and $(w_{0},k_{0})$ correspond under the bijection $(X\ti \Z) \rl (\W\ti \Z)$. 
The equivalence class of $(w_{0},k_{0})$ is the set $E$ from Corollary \ref{vgs} (a).  
If  $(\gl,  m\ell)\sim (\mu,  mr)$, then $(\mu,  mr)$  
corresponds to  a pair $(w_{i},k_{i})$ such that $m$ divides $k_i -k_0$.  Hence $(w_{0},k_{0})\ap (w_{i},k_{i})$ by Corollary \ref{vgs} (d) and the result follows. 
\epf
 \bexa \label{E3}{\rm  
 Let $(n,m) =(2,3)$.  The Hasse diagrams for $[X\ti m\Z]_\ap$  and $[X\ti \Z]$ 
are isomorphic by Theorem \ref{vss}.  
We give the portion of the Hasse diagram of  $[X\ti 3\Z] $ in degrees 0-6.
  The rotation number $k$ is given as a superscript and is omitted if $k=0$.   The degree, as defined in \eqref{e67n}  is given in the last row. Because of the periodic nature of this poset, the part in degree $i+6$ is obtained  from the degree $i$ part by increasing the rotation number by 6. So the portion below determines the entire graph.   
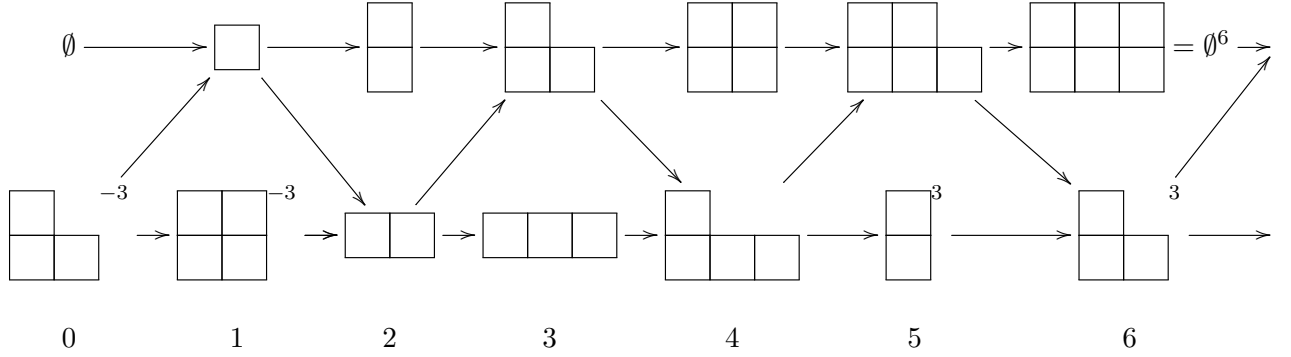
\begin{equation*}
\xymatrix@C=1pc@R=1pc{\emptyset\ar@{->}[r] &
\ydiagram{1}\ar@{->}[ddr]&
\ydiagram{1,1} \ar@{<-}[l] \ar@{->}[r]& \ydiagram{0+1,2}\ar@{->}[r] \ar@{->}[ddr]   &
\ydiagram{2,2} \ar@{->}[r]&\ydiagram{0+2,3} \ar@{->}[r] \ar@{->}[ddr]& \ydiagram{3,3}=\emptyset^{6}\ar@{->}[r]&
\\ \\
 \ydiagram{0+1,2}^{-3}\ar@{->}[r] \ar@{->}[uur]&\ydiagram{2,2}^{-3}\ar@{->}[r]&
\ydiagram{2}\ar@{<-}[l] \ar@{->}[r] \ar@{->}[uur]&\ydiagram{3}\ar@{->}[r]&
\ydiagram{0+1,3}\ar@{->}[r] \ar@{->}[uur]&\ydiagram{1,1}^3\ar@{->}[r]& \ydiagram{0+1,2}^3\ar@{->}[r]\ar@{->}[uur]&\\
0&1&2&3&4&5&6}
\end{equation*}
 The most interesting behavior occurs when the rotation number changes after a morphism is applied.  Here are some examples 
By  \eqref{1no} 
we obtain 
\be \label{sb1}\F(\gr_{\gep_2 -\gd_1})[(3,i),0)]=[t_{\gep_1 -\gd_1}(i,0),3] =[((i,1),{3})]=[(i,1)^{3}] \mbox{ for } i=1,2.\ee
Together with the more obvious morphism $(1,1)^3 \lra (2,1)^3$, this accounts for the extra arrows in the lower right that do not appear in the Hasse diagram for the poset $X$.
Using   \eqref{2no} 
$$\F(\gr_{\gd_3-\gep_1})[(i,0),0)]=\F(\gr_{\gd_3-\gep_1})[(3,i),{-3})]=[t_{\gd_3-\gep_2 }(3,i),{-3}]  =[(2,i)^{-3}]  \mbox{ for } i=1,2.$$
The diagram shows only positive morphisms, but the last equation is equivalent to
$$\F(\gr_{\gep_1-\gd_3})[(2,i)^{-3}] =[(i,0),0)] \mbox{ for } i=1,2.$$
This accounts for the  extra arrows in the lower left  of the diagram.  In general the proof of Lemma \ref{vgt} shows that the arrows that do not appear in the Hasse diagram for  $X$ arise from pseudocorners.
}\eexa

\section{Borel subalgebras} \label{akmls} 

\noi In this Section $\fgo$ will denote the Lie superalgebra $\fsl(n|m)$. 
Set $r=m+n-1$. 
Let $\fho$ be the subalgebra of diagonal matrices in  $\fgo$. 
Suppose that $ A^0  = (a^0_{ij})$
is the Cartan matrix of $\fgo$ corresponding to the shuffle $\gs$.
\\ \\
 {\bf Notation: }We assume $w, \gl$ and $\gs$ are related as in \eqref{bji}.

\subsubsection{Outline} \label{spafn}   
Our goal is to construct Borel subalgebras as in \eqref{f2}.    
Let $\pigo$ be the distinguished set of simple roots for $\fgo$, $\fbo$ the distinguished Borel subalgebra   
and let  $\gb=\gep_n-\gd_1$ be the unique odd root in $\pigo$.  
For $\gl\in X$ we have a corresponding Borel $\fbo(\gs)$ as in \eqref{bji} with simple roots 
$\pigo(\gs) $
and  highest root $\ghs$.
Set $\ga_0(\gs)=\ogd -\ghs$ where $\ogd$ is the basic imaginary root, see \eqref{189}.  In the affinization $\Lgh$ set 
$ \Pi(\gs) =\pigo(\gs)\cup\{\ga_0(\gs)\}.$  The corresponding extended Dynkin-Kac diagram is obtained from the Dynkin-Kac diagram for $\fgo$
arising from the simple roots $\pigo(\gs) $ 
by adjoining a node corresponding to the root $\ga_0(\gs)$.  Under the conditions of \eqref{f2} 
by deleting a node of the extended diagram which is adjacent to $\ga_0(\gs)$ 
 we obtain a set of simple roots for a
different subalgebra $\fgo'$ of $\Lgh$ which is isomorphic to $\fgo$.   This allows us to perform further odd reflections that do not arise from $\fgo$.  We can  repeat this process to obtain  subalgebras $\fgo(k)$ for $k\in \Z$.   
Below we provide further details of this procedure. In particular we need generators for the algebras to explain how $\fgo'$ is embedded in  $\Lgh$. 

\subsubsection{Definition of Kac-Moody superalgebras} \label{ss5.1}
For convenience we recall the definition of  Kac-Moody Lie superalgebras (also known as a contragredient Lie superalgebras). Let $A = (a_{ij})$ be an $N \times N$ symmetric matrix of rank $l$ with
entries in the base-field $\ttk$ and let  $(\mathfrak{h}, \Pi, \Pi^{\vee}),$ be a  minimal realization of $A$.   By definition   $\mathfrak{h}$ is a vector space of dimension $2N-\ell$ and \be
\label{pipi} \Pi = \{\alpha_1, \ldots , \alpha_N\}, \quad \Pi^{\vee}
= \{h_1, \ldots, h_N\}\ee are linearly independent subsets of
$\mathfrak{h}^*, \mathfrak{h}$ respectively, such  that \be \label{re1}
\alpha_j(h_i) =a_{ij}.\ee
If  $\tau$ is a  subset of $\{1,2, \ldots, N\}$, the {\it half-baked Lie
superalgebra} $\widetilde{\fg} = \widetilde{\mathfrak{g}}(A,\tau)$ associated to $(A,\tau)$ is generated by the
vector space $\mathfrak{h}$ and  elements $e_i, f_i, h_i \; (i = 1, \ldots,
N)$ with $h_i \in  \fh$ subject to the relations
\begin{equation}
\label{re3}
   [ e_i, f_j ]  =\delta_{ij} h_i,
\ee
\be\label{re4}
   [ h, h'] = 0,
\ee
\be\label{conLS}
   [ h, e_i] =       \alpha_i (h) e_i, \quad
    [h, f_i]  =    - \alpha_i (h) f_i
 \ee
 for $i, j = 1, \ldots , N$ and $h,h' \in \mathfrak{h}.$
Note that Equations (\ref{re1}) and (\ref{conLS}) imply that \be
\label{rels}
   [ h_i,e_j]       = a_{ij}e_j, \quad
    [h_i, f_j]      = - a_{ij}f_j.
\ee
 The $\mathbb{Z}_2$-grading on $\widetilde{\mathfrak{g}}(A,\tau)$ is given
by
\[ \begin{array}{lll}
\deg h_i = 0, & \deg e_i = \deg f_i = 0 & \mbox{for} \;  i \notin
\tau \\
              & \deg e_i = \deg f_i = 1 & \mbox{for} \;  i \in \tau.
\end{array} \]
There is a unique ideal $\mathfrak{r}$ of
 $\widetilde{\mathfrak{g}}$ which is maximal among all ideals of
 $\widetilde{\fg}$ which intersect $\mathfrak{h}$ trivially and we set \be \label{clsdef} \mathfrak{g}(A,\tau) =
\widetilde{\mathfrak{g}}(A,\tau)/\mathfrak{r}.\ee
\subsection{$\fgo$ as a KM Lie superalgebra} \label{afn1}
\noindent The algebra  $\fgo$ can be constructed as a  KM  Lie superalgebra   $\fgo = \fg(A^0, \tau^0)$. 
We briefly recall some details. Let  ${\stackrel{{\rm o}}{\fb}}{}(\gs)$ be the Borel subalgebra of $\fgo$ corresponding to the shuffle $\gs$. We are given subsets 
\be \label{pihid}\pigo = \{\alpha_1, \ldots, \alpha_r\} \;\; \mbox{and} \;\;
\stackrel{\rm o}{\Pi^\vee} = \{H_1,\ldots, H_r\}.\ee
of $\fho$ and $\fho^*$ respectively.  The set 
$\pigo$ is the set of simple roots of ${\stackrel{{\rm o}}{\fb}}{}(\gs)$  and the $\ga_i$ are given as in \eqref{pgs} by $\ga_i=\gep_{\gs(i)}-\gep_{\gs(i+1)}$. For $i\in [r]$ let   $E_i, F_i$ be root vectors with weights $\ga_i,-\ga_i$  respectively and set    $H_i=[E_i, F_i]$.  Then $\fgo$ is generated by the $E_i, F_i$ and the Cartan subalgebra $\fho$ of diagonal matrices. These elements satisfy the relations
(18.2.2)-(18.2.4) from \cite{M101}.
The subalgebra $\fb(\gs)$ is generated by the $E_i $ and  $\fho.$  
The  matrix $A^0$ has entries  $a^0_{ij}=\alpha_i(H_j)$.  
Also all elements of $\fho$ are even, and for $ i \in [r]$, $E_i$ and $F_i$ are odd if and only if $i \in \tau^0$. By  \cite{M101} Theorem 5.3.5 we have

\bl \label{pih} A minimal realization of  $A^0$ is given by   
 $(\fho, \pigo, \pigov)$.  We have $\fgo \cong \fg(A^0, \tau^0)$.
\el 
\noi 
If we perform an odd reflection using an odd simple root in $\pigo$ we obtain a new set of generators for $\fgo$ and a new Cartan matrix $A_1^0$.  The relationship between $A^0$ and $A^0_1$ is explained in \cite{M101} Section 3.5 and  Lemma 5.3.4.
\subsection{Affinization} \label{afn}
Starting from $A^0$ we construct
another matrix $A = (a_{ij})$ with an extra row and column labeled
by 0. 
 There is a unique vector $v$ such that the augmented matrix below $A$ has all row sums and column sums equal to zero.  
The matrix $A$ can also be decomposed using 
a new matrix  $\fag$ and a suitable vector $w$.  Thus 
\be \label{mtr} A=\begin{tabular}{|c||ccc|} \hline
0&  &$v^\trans$&
\\ \hline\hline &     &                          & \\ $v$
&&$A^0$
 & \\
 & &                              & \\
  \hline
\end{tabular}
=\begin{tabular}{|ccc||c|} \hline
&  &&\\
 &    $\fag$  &                          & $w$\\ 
&& & \\
 \hline\hline & $w^\trans$&                              & 0\\
  \hline
\end{tabular}\ee 
It is clear that each of $A^0,  \fag$ determines the other, 
Let $\tau$ be a subset of $\{ 0,\ldots, r \}$,  and set $ \tau^0  =\gt \cap \{ 1,\ldots, r \}$,  
$\fto = \tau \cap \{ 0,\ldots, r-1 \}.$ We give an
 explicit construction of the  KM  Lie superalgebra
 $\fg(A,\tau)$.\\
\\
 The {\it loop algebra} $\Lgo$ is the Lie superalgebra with graded
 components, $(i = 0,1)$
 \[ \Lgo_i = \ttk[t,t^{-1}] \otimes \fgo_i.\]
Next the 
invariant bilinear form 
 $( \; , \;)$ 
on $\fgo$ is used to form a
 one dimensional central extension $\Lgo \oplus \ttk c $ of $\Lgo$, \cite{M101} Lemma 18.2.1. The product is given by 
\be \label{mjr} [t^m \otimes x, t^n \otimes y] = t^{m+n} \otimes [x,y] + m(x,y)
\delta_{m,-n} c. \ee
Define a derivation $d$ on $\Lgo \oplus \ttk c$ such that $d(c)=0$ and $d$ acts as 
$t\frac{d}{dt}$ on the loop algebra. 
The {\it affinization} of $\fgo$ is 
 \[ \Lgh = \Lgo \oplus \ttk c \oplus \ttk d. \]
 We make $\Lgh$ into a Lie superalgebra with $d$ even, such
 that the product on $\Lgh$ extends the product already defined on
$\Lgo \oplus \ttk c,$ by requiring that
$ [d,x] = d(x)$
for all $x \in \Lgo \op \ttk c$.
Let  $\theta =
 \sum^r_{i=1}\alpha_i$  be the highest root of $\fgo$ and choose elements 
$F_0 \in \; \fgt, E_0 \in \;
\fgm,$ suitably normalized as in \cite{M101} Equation (18.2.12) and such  that $(E_0 ,F_0)=1$.  
Set $$\fh = (1 \otimes \fho) \op \ttk c \oplus \ttk d,$$ and extend
the elements $\alpha_1, \ldots, \alpha_r \in \; \fhs$ to linear
forms on $\fh$ by setting $\alpha_i(c) = \alpha_i(d) = 0$ for $ i \in  [r]$. Define $\overline{\delta} \in \fh^*$ by
\be\label{189} \overline{\delta}(c) = 0, \;\; \overline{\delta}(d) = 1, \;\; \overline{\delta}(h) = 0 \;\; \mbox{for} \;\; h
\in \; \fho . \ee
In  $\Lgh$, set for $i\in [r]$
\[ e_i = 1 \otimes E_i, \quad f_i = 1 \otimes F_i, \quad h_i = 1 \otimes H_i. \]
Set $e_0 = t \otimes E_0, f_0
= t^{-1} \otimes F_0.$ Then by \cite{M101} Lemma 8.3.2, 
\be\label{efo}
 h_0 = [e_0, f_0]  =  [E_0, F_0] + (E_0, F_0)c \nonumber \\
           =  -(1 \otimes h_\theta)+ c.
\ee
We set $\alpha_0 = \overline{\delta} - \theta \in \; \fh^*$,  and $h_i = 1 \otimes H_i \in \fh$ for $i \in [r]$.
 Then set 
\be\label{187}
\Pi =
 \{\alpha_0, \alpha_1, \ldots, \alpha_r\}, \quad\mbox{ and } \quad
 \Pi^\vee =
 \{h_0, h_1, \ldots, h_r\}.\ee  The triple $(\fh, \Pi, \Pi^\vee)$ is a minimal realization of $A$ and  $\Lgh \cong \fg(A,\tau)$,
\cite{M101}, Lemma 18.2.2 and we have by  \cite{M101} Theorem 18.2.5. 
  \bt \label{185} With the above notation, $\Lgh \cong \fg(A,\tau)$.  This algebra is generated by $e_0, e_1,
\ldots, e_r,$ $ f_0, f_1, \ldots, f_r$ and $\fh$.
\et \noi 
\br\label{186} {\rm 
The proof of Theorem \ref{185} shows more. The Cartan subalgebra $\fh$ appears in the definition of both $\Lgh$ and $\fg(A,\tau)$ and  the isomorphism is the identity on $\fh$. Also  the generators $e_0, e_1, \ldots, e_r$ and $f_0, f_1, \ldots, f_r$  of $\Lgh$ correspond under the isomorphism to the canonical generators of $\fg(A,\tau)$  as a  KM  Lie superalgebra. 
}\er

\subsection{A fundamental construction} \label{afc}
Let 
\be\label{efx}
\pou= \{\alpha_0, \ldots, \alpha_{r-1}\}, \; 
\pvu =
 \{h_0, \ldots, h_{r-1}\} \; \mbox{and} \; {\underline{\fh}^0} =\span \pvu .\ee
Then the triple $(\underline{\fh}^0, \pou, \pvu)$ is a minimal realization of $\fag$.  If $
\fgu = \fg(\fag, \fto)$, then  by  the proof of Theorem \ref{185} we have ${\widehat{L}(\fgu})\cong \fg(A,\tau)$.  
Moreover viewed as a subalgebra of $\fg(A,\tau)$ by means of this isomorphism, $\fgu$ is generated by $ \{e_0, \ldots, e_{r-1}\},  \{f_0, \ldots, f_{r-1}\}$ and $\underline{\fh}^0$, compare Remark \ref{186}.
Then  $\fgu$ is  the subalgebra of $\Lgh$ whose 
Dynkin-Kac diagram is obtained from the diagram 
for $\Lgh$ by deleting the node corresponding to the simple root $\gep_{\gs(r)}-\gep_{\gs(r+1)}$ as explained in Subsection \ref{spafn}.  
Then $\fgu\cong \fgo$. 
\noi We say that the Borel subalgebra of $\Lgh$ with simple roots $\Pi$ as in \eqref{187} {\it is the extension to} $\Lgh$  of the Borel subalgebra of ${\stackrel{{\rm o}}{\fg}}$ with simple roots ${\stackrel{{\rm o}}{\Pi}}$ as in \eqref{pihid}.  Note that the {\it extending root} $\alpha_0 = \overline{\delta} - \theta$ is determined by ${\stackrel{{\rm o}}{\Pi}}.$  In a similar way 
$\fb$ is the extension to $\Lgh$  of the Borel subalgebra of $\fgu$ with simple roots $\pou$.
\\ \\
Recall that for $\fgo=\fgo(0)$ and a given shuffle $\gs$, the set of simple roots ${\stackrel{{\rm o}}{\Pi}}(\gs)$ is given by  
\eqref{pgs}. For Borels and simple roots in $\fgo(k) $ we use similar notation reflecting the fact that $\fgo(k) \cong \fgo(0)$. We call the resulting expressions the {\it local names } for these objects.   
However if we want to express the simple roots of ${\stackrel{{\rm o}}{\fb}}{}(\gs,k)$ in terms of the distinguished roots and the extending root  of $\fgo(0)$, we use {\it global names}. For our purposes the latter are less useful than local names.  Moreover they are harder to compute, involving powers of $\nu, \eta$ and multiples of $\ogd$.  Therefore we restrict our use of global names to one example, Example \ref{glh}.

\subsection{Definition  of the algebras from \eqref{f2}}  \label{afz}
By repeating the construction of the previous Subsection, we obtain  various subalgebras of $\Lgh$ which are all isomorphic to $\fgo$. It will no longer be necessary to consider explicit generators for these subalgebras.   Instead we work with sets of simple roots.  It may be worth noting that if a root is not imaginary then the corresponding root space has dimension one, \cite{M101} Corollary 18.2.4.  Also an imaginary root never occurs as a simple root.  \\ \\  
Starting from ${\stackrel{{\rm o}}{\fb}}{}(\gl,k)$
 and  ${\stackrel{\rm o}{\fg}}(k)$, we explain how to construct the algebras from \eqref{f2}.
We consider the two cases where the last letter of a word is moved to the front, cases (a) and (d) of Lemma \ref{rogs}.
If $\gl \in X(\rda)$,  assume 
 ${\stackrel{{\rm o}}{\fb}}{}(\gl,k)$  is the Borel subalgebra  of
 $\fg(A^0, \gt^0)= {\stackrel{\rm o}{\fg}}(k)$ with simple roots $\pigo$ as in  \eqref{pihid}.  Then 
 ${\stackrel{{\rm o}}{\fb}}{}(\gl^\mr,k+m)$ is the  Borel subalgebra of $\fg(\fag, \fto) ={\stackrel{\rm o}{\fg}}(k+m)$
with simple roots $\pou$ as in \eqref{efx}.
\\ \\
Now suppose $\gl \in X(\cpa)$ and   assume that 
 ${\stackrel{{\rm o}}{\fb}}{}(\gl,k)$  is the Borel subalgebra  of
 $\fg(\fag, \fto)= {\stackrel{\rm o}{\fg}}(k)$ with simple roots  $\pou$ as in \eqref{efx}.
Then 
 ${\stackrel{{\rm o}}{\fb}}{}(\gl^\opc,k-n)$ is the  Borel subalgebra of $\fg(A^0, \gt^0) ={\stackrel{\rm o}{\fg}}(k-n)$
with simple roots $\pigo$ as in \eqref{pihid}.    It remains to define the corresponding algebras from \eqref{f2} when 
$\gl \in X(\rpa)$ or $\gl \in X(\cda)$, but it is clear that this should be done by reversing the procedures given above.  We leave the details to the reader. 
\\  \\Since $n,m $ are coprime, we can 
start from $k=0$ and repeating the above constructions obtain subalgebras  ${\stackrel{{\rm o}}{\fb}}{}(\gl,k)$  of ${\stackrel{\rm o}{\fg}}(k)$ for all $k\in\Z$.  
At first this is done only for $\gl$ satisfying the  conditions of \eqref{f2}, but by using odd reflections we have them for all $\gl\in X$.

\subsection{The equivariant bijection }  \label{bfz}

Let $\cB$ be the set of Borel subalgebras of $\Lgh$ that are extensions  of Borels in the subalgebras $\fgo(k)$ for  $k\in \Z.$  

\bp \label{crb} There is well defined an equivariant bijection $x:{[\stackrel{{\rm o}}
{\cB}} \ti \Z] \lra \cB$ where $[{\stackrel{{\rm o}}{\fb}}{}(\gs),k]$  maps to the Borel subalgebra of $\fg$ extended from ${\stackrel{{\rm o}}{\fb}}{}(\gs,k)$.
\ep
\bpf Since the equivalence relation $\sim$ is generated by \eqref{1kr}  and \eqref{2kr},  to show $x$ is well-defined, it is enough to check the following 
\bi \itema  
If $\gl \in X(\rda)$ the Borel subalgebras ${\stackrel{{\rm o}}{\fb}}{}(\gl^\mr,k+m)$ and ${\stackrel{{\rm o}}{\fb}}{}(\gl,k)$
 extend to the same Borel subalgebra of $\Lgh$.
\itemb  
If  $\gl \in X(\cpa)$ the Borel subalgebras ${\stackrel{{\rm o}}{\fb}}{}(\gl^\opc,k-n)$ and ${\stackrel{{\rm o}}{\fb}}{}(\gl,k)$
 extend to the same Borel subalgebra of $\Lgh$.
\ei  This is straightforward.
Note that equivariance follows from Lemma \ref{ymn} and the remarks earlier in this Subsection. 
The map ${x:[\stackrel{{\rm o}}{\cB}} \ti \Z] \lra \cB$ is onto by definition of $\cB$.  To show it is injective, suppose
$x[\fbo(w,k)] =\fb$.  Then 
DK diagram $D$ for $\fb$ determines a set of $m+n$ words $w_i$ as in Subsection \ref{dta}, see  also Example \ref{iob}.
Any Borel  that extends to $\fb$ has a diagram obtaained by deleting a node from $D$. But we have shown in Corollary \ref{vgs} that the equivalence class of $\fbo(w_0,k_0)$ has size $m+n$ and consists of the set $\{\fbo(w_i,k_i)| 0 \le  i < m+n\}.$  
Therefore the preimage of $\fb$ under $x$ consists of this single equivalence class.
\epf
\bt \label{srb} The groupoid  $
\mathfrak  T_{iso}$ acts on $\cB$ by odd reflections and we have a commutative 
diagram 
\[
\xymatrix@C=2pc@R=1pc{
\mathfrak  T_{iso}\ar@{<->}_= [dd]&&
\cS(\cB) \ar@{<-}[ll]: \ar@{<->}^\cong[dd]&\\ \\
\mathfrak  T_{iso} \ar@{-}[rr] &&
\cS[{\stackrel{{\rm o}}{\cB}} \ti \Z]\ar@{<-}[ll] &}\]
\et
\bpf This  follows from Proposition \ref{crb}.  
\epf \noi 

\bexa \label{glh}{\rm In Example \ref{glk} we carried out some computations with partitions (equivalently with words), indicating how to pass from  ${ \boldsymbol \gl}$ to 
$\emptyset$  by deleting a row and a column and using a single odd reflection.  For convenience we recall some  details when $n=3,  m=4$.  Thus    ${\gl= \boldsymbol \gl}=\ttr^{}  \ttd^{2}\ttr^{3}\ttd^{}$, $\gl^\mr=\ttd^{}\ttr^{}  \ttd^{2}\ttr^{3}$ and  
$\mu=\ttr^{}  \ttd^{3}\ttr^{3}.$ 
Now by \eqref{1kr} 
$(\gl,0)\sim (\gl^\mr, 4)$. Then using an odd reflection, we have 
$r_\ga(\gl^\mr,4)=(\mu, 4)$ 
and  by \eqref{2kr} 
$(\mu,4)\sim (\mu^\omc, 7)= (\emptyset, 7)$.
In the table we give the global names for the simple roots of the Borel subalgebras corresponding to the above calculation.    
}\eexa

\[ \begin{tabular}{|c||c|} \hline
Borel  Subalgebra &\;\; Global Names for Simple Roots \;\;  \\ \hline
$\fbo(\gl,0) $
&
$\gd_1-\gep_1,\;\; \gep_1-\gep_2,\;\;  \gep_2-
\gd_2,\;\; \gd_2 - \gd_3,\;\; \gd_3 -\gd_4,\;\; \gd_4-\gep_3.$
\\ \hline
$\fbo(\gl^\mr, 4)$ 
&
$\ogd - \gd_1+\gep_3,\;\;
\gd_1-\gep_1,\;\; \gep_1-\gep_2,\;\;  \gep_2-
\gd_2,\;\; \gd_2 - \gd_3,\;\; \gd_3 -\gd_4. 
$
\\ \hline
$\fbo(\mu,\;4)$ &  $-\ogd +\gd_1-\gep_3,\;\;
\ogd +\gep_3-\gep_1,\;\; \gep_1-\gep_2,\;\;  \gep_2-
\gd_2,\;\; \gd_2 - \gd_3,\;\; \gd_3 -\gd_4. 
$
\\ \hline
$\fbo(\emptyset, 7)$  & $
\ogd +\gep_3-\gep_1,\;\; \gep_1-\gep_2,\;\;  \gep_2-
\gd_2,\;\; \gd_2 - \gd_3,\;\; \gd_3 -\gd_4,\;\; \ogd +\gd_4- \gd_1.
$\\ \hline
\end{tabular}\]
The first list of roots is given by \eqref{pgs}.  To pass from the first/third to the second/last row use the affinization procedure.  The third list is obtained by odd reflection using the root 
$\ogd - \gd_1+\gep_3.$


\bexa \label{E4}{\rm We consider certain Borel subalgebras of $\Lgh$, when $\fgo=\fsl(2|3)$. 
Let $\gs, \gs_1$ be the shuffles with
$${ \boldsymbol \gs}
= ( 1',1, 2',3', 2 ) \mbox{ and } { \boldsymbol \gs}_1 
= ( 1',2', 1, 3', 2).$$
Let $\fb(\gs)$ and $\fb(\gs_1)$ be the Borel subalgebras of ${\stackrel{{\rm o}}{\fg}}$ corresponding to $\gs$ and $ \gs_1$.  
Note that $\fb(\gs)$ and $\fb(\gs_1)$  are related by odd reflections using the roots
$\pm(\gep_1-\gd_2).$ We give two examples of the procedure defined above.  In the first
(resp. second) the Dynkin-Kac diagram for $\fb(\gs)$ (resp.  $\fb(\gs_1)$) is given on the left, the corresponding diagram for the affinization is in the middle and if $\gt$ (resp. $\gt_1$) is the shuffle obtained by the procedure, the Dynkin-Kac diagram for $\fb(\gt)$ (resp.  $\fb(\gt_1)$) is given on the  right.    
Note that $\fb(\gt)$ and $\fb(\gt_1)$  are related by odd reflections using the roots
$\pm(\gep_2-\gd_2).$

\Bc
\setlength{\unitlength}{0.9cm}
\begin{picture}(10,1)(-0.5,1.70)
\thinlines
  \linethickness{.09mm}
\put(-2.4,1.81){\line(1,0){0.74}}
\put(-1.4,1.81){\line(1,0){0.74}}
\put(-0.4,1.81){\line(1,0){0.74}}

\put(-3.3,1.3){$\scriptstyle 1'$}
\put(-2.3,1.3){$\scriptstyle 1$}
\put(-1.3,1.3){$\scriptstyle 2'$}%
\put(-0.3,1.3){$\scriptstyle 3'$}
\put(.73,1.3){$\scriptstyle 2$}\put(0.3,1.7){{$\otimes$}}

\put(-2.7,1.7){{$\otimes$}}
\put(-1.7,1.7){{$\otimes$}}\put(-0.59,1.81){\circle{0.27}}
\put(9,1.81){\line(1,0){0.74}}
\put(10,1.81){\line(1,0){0.74}}
\put(11,1.81){\line(1,0){0.74}}

\put(8.1,1.3){$\scriptstyle 1$}
\put(9.1,1.3){$\scriptstyle 1'$}
\put(10.1,1.3){$\scriptstyle 2$}
\put(11.27,1.3){$\scriptstyle 2'$}\put(12.13,1.3){$\scriptstyle 3'$}%

\put(8.7,1.7){{$\otimes$}}
\put(9.7,1.7){{$\otimes$}}
\put(10.7,1.7){{$\otimes$}}\put(11.87,1.81){\circle{0.27}}%
\put(3.7,0.82){{$\otimes$}}
\put(4.87,0.93){\circle{0.27}}
\put(4.2,2.58){$\otimes$}
\put(2.95,1.7){{$\otimes$}} 
\put(5.46,1.7){$\otimes$}
\put(5.53,1.94){\line(-5.2,3){1.09}}
\put(3.19,1.92){\line(5,3.2){1.06}}
\put(3.75,1.01){\line(-4,5){0.54}}
\put(4.95,1.01){\line(4,5){0.54}}
\put(3.3,1.07){$\scriptstyle 1$}
\put(4.3,.42){$\scriptstyle 2'$}
\put(5.3,1.07){$\scriptstyle 3'$}
\put(3.17,2.28){$\scriptstyle 1'$}
\put(5.27,2.28){$\scriptstyle 2$}
\put(4.0,.93){\line(1,0){0.74}}
\end{picture}
\Ec 
\vspace{1cm}

\Bc
\setlength{\unitlength}{0.9cm}
\begin{picture}(10,1)(-0.5,1.70)
\thinlines
  \linethickness{.09mm}
\put(-2.47,1.81){\circle{0.27}}
\put(-1.7,1.7){{$\otimes$}}
\put(-.62,1.7){$\otimes$}\put(.36,1.7){$\otimes$}
\put(-2.4,1.81){\line(1,0){0.74}}
\put(-1.4,1.81){\line(1,0){0.74}}
\put(-0.4,1.81){\line(1,0){0.74}}
\put(-3.1,1.3){$\scriptstyle 1'$}
\put(-2.1,1.3){$\scriptstyle 2'$}
\put(-1.1,1.3){$\scriptstyle 1$}
\put(.1,1.3){$\scriptstyle 3'$}
\put(1.1,1.3){$\scriptstyle 2$}
\put(3.7,0.82){{$\otimes$}}
\put(4.7,0.82){{$\otimes$}}
\put(4.2,2.58){$\otimes$}
\put(3.07,1.87){\circle{0.27}} 
\put(5.46,1.64){$\otimes$}
\put(5.53,1.94){\line(-5.2,3){1.09}}
\put(3.19,1.92){\line(5,3.2){1.06}}
\put(3.75,1.01){\line(-4,5){0.54}}
\put(4.95,1.01){\line(4,5){0.54}}
\put(3.3,1.07){$\scriptstyle 2'$}
\put(4.3,.42){$\scriptstyle 1$}
\put(5.3,1.07){$\scriptstyle 3'$}
\put(3.17,2.28){$\scriptstyle 1'$}
\put(5.27,2.28){$\scriptstyle 2$}
\put(4.0,.93){\line(1,0){0.74}}

\put(9,1.81){\line(1,0){0.74}}
\put(10,1.81){\line(1,0){0.74}}
\put(11,1.81){\line(1,0){0.74}}

\put(8.1,1.3){$\scriptstyle 1$}
\put(9.1,1.3){$\scriptstyle 1'$}
\put(10.1,1.3){$\scriptstyle 2'$}
\put(11.27,1.3){$\scriptstyle 2$}
\put(12.13,1.3){$\scriptstyle 3'$}%

\put(8.7,1.7){{$\otimes$}}
\put(9.87,1.81){\circle{0.27}}
\put(10.7,1.7){{$\otimes$}}
\put(11.76,1.7){$\otimes$}
\end{picture}
\Ec 
\vspace{1cm}
}\eexa



\end{document}